\documentclass{article}

\usepackage{amscd}      
\usepackage{amssymb}     
\usepackage{amsmath}     
\usepackage{xypic}      
\LaTeXdiagrams          
\usepackage[all,v2]{xy}      
\xyoption{2cell}
\UseAllTwocells
\xyoption{frame}
\CompileMatrices
\allowdisplaybreaks[3] 
\usepackage{theorem}

\newtheorem{prop}{Proposition}[section]  
\newtheorem{lem}[prop]{Lemma}

\newtheorem{cor}[prop]{Corollary}
\newtheorem{them}[prop]{Theorem}

\newtheorem{caveat}{Caveat}

\theorembodyfont{\upshape}

\newtheorem{defn}[prop]{Definition}

\newtheorem{numrmk}[prop]{Remark}

\newtheorem{numex}[prop]{Example}

\newtheorem{example}{Example}

\newtheorem{rmk}{Remark}

\newenvironment{pf}{\begin{trivlist}\item[]{\sc Proof.}}%
            {\nolinebreak $\Box$ \end{trivlist}}
\newcommand{\unsure}[1]{{\footnotesize #1}}

\newcommand{\noprint}[1]{}

\newcommand{\comment}[1]{{\marginpar{\footnotesize #1}}}

\renewcommand{\tilde}{\widetilde}

\newcommand{\toto}{\rightrightarrows}
\newcommand{\qed}{{\nolinebreak $\,\,\Box$}}

\newcommand{\an}{{\mbox{\tiny an}}}
\newcommand{\bigsite}{{\mbox{\tiny big}}}

\newcommand{\upst}{^{\ast}}

\newcommand{\lst}{_{\ast}}

\newcommand{\com}{^{\scriptscriptstyle\bullet}}
\newcommand{\lcom}{_{\scriptscriptstyle\bullet}}

\newcommand{\X}{\mathop{\sf X}}

\newcommand{\GG}{{\mathfrak G}}

\newcommand{\Gg}{{\mathfrak g}}

\newcommand{\LL}{{\mathfrak L}}

\newcommand{\hh}{{\mathbb H}}
\newcommand{\kk}{{\mathbb K}}

\newcommand{\cc}{{\mathbb C}}
\newcommand{\rr}{{\mathbb R}}

\renewcommand{\O}{{\cal O}}

\newcommand{\gG}{{\cal G}}
\newcommand{\fF}{{\cal F}}

\newcommand{\del}{\partial}

\newcommand{\resto}{{\,|\,}}

\newcommand{\rank}{\mathop{\rm rank}\nolimits}

\newcommand{\id}{\mathop{\rm id}\nolimits}

\newcommand{\comp}{\mathbin{{\scriptstyle\circ}}}

\newcommand{\ldiag}[1]%
       {\makebox[0cm]{${\scriptstyle#1}\downarrow\phantom{\scriptstyle#1}$}}
\newcommand{\ldiagup}[1]%
       {\makebox[0cm]{${\scriptstyle#1}\uparrow\phantom{\scriptstyle#1}$}}
\newcommand{\rdiag}[1]%
       {\makebox[0cm]{$\phantom{\scriptstyle#1}\downarrow{\scriptstyle#1}$}}
\newcommand{\sediagr}[1]%
       {\makebox[0cm]{$\phantom{\scriptstyle#1}\searrow{\scriptstyle#1}$}}
\newcommand{\nediagr}[1]%
       {\makebox[0cm]{$\phantom{\scriptstyle#1}\nearrow{\scriptstyle#1}$}}
\newcommand{\rdiagup}[1]%
       {\makebox[0cm]{$\phantom{\scriptstyle#1}\uparrow{\scriptstyle#1}$}}
\newcommand{\swdiag}[1]%
       {\makebox[0cm]{$\phantom{\scriptstyle#1}\swarrow{\scriptstyle#1}$}}
\newcommand{\sediag}[1]%
       {\makebox[0cm]{${\scriptstyle#1}\searrow\phantom{\scriptstyle#1}$}}
\newcommand{\nediag}[1]%
       {\makebox[0cm]{${\scriptstyle#1}\nearrow\phantom{\scriptstyle#1}$}}

\newcommand{\longiso}{\stackrel{\textstyle\sim}{\longrightarrow}}
\newcommand{\iso}{\stackrel{\sim}{\rightarrow}}

\newcommand{\doublearrowstack}[2]%
                      {{{{\scriptstyle#1}\atop{\textstyle\longrightarrow}}\atop{{\textstyle\longrightarrow}\atop{\scriptstyle#2}}}}
\newcommand{\rightleftarrowstack}[2]%
                      {{{{\scriptstyle#1}\atop{\textstyle\longrightarrow}}\atop{{\textstyle\longleftarrow}\atop{\scriptstyle#2}}}}
\newcommand{\leftrightarrowstack}[2]%
                      {{{{\scriptstyle#1}\atop{\textstyle\longleftarrow}}\atop{{\textstyle\longrightarrow}\atop{\scriptstyle#2}}}}

\newcommand{\overtoparrow}%
{\makebox[0cm]{\beginpicture
\setcoordinatesystem units <.8cm,.4cm> point at 0 0
\setplotarea x from -3 to 3, y from 0 to 1
\setquadratic
\plot -3 0 0 1 3 0 /
\put{\vector(3,-1){0}}[Bl] at 3 0
\endpicture}}

\newcommand{\underbottomarrow}%
{\makebox[0cm]{\beginpicture
\setcoordinatesystem units <.8cm,.4cm> point at 0 0
\setplotarea x from -3 to 3, y from 0 to 1
\setquadratic
\plot -3 1 0 0 3 1 /
\put{\vector(3,1){0}}[Bl] at 3 1
\endpicture}}

\newcommand{\ses}[5]%
{0\longrightarrow#1\stackrel{#2}{ \longrightarrow}#3\stackrel{#4}{
\longrightarrow}#5\longrightarrow0}

\newcommand{\dt}[6]%
{#1\stackrel{#2}{longrightarrow}#3 \stackrel{#4}{\longrightarrow}#5
\stackrel{#6}{\longrightarrow} #1[1]}  
 
\newcommand{\cat}[1]%
{(\mbox{\rm #1})}

\def\Label#1{\label{#1}{\tt [#1]}\phantom{h}}
\def\Label{\label}
\def\comment{\noprint}
\def\unsure{\noprint}

\title{On the de Rham Cohomology of Differential and Algebraic Stacks}
\author{Kai A. Behrend}

\begin{document}
\sloppy
\maketitle
\begin{abstract}
We introduce the notion of {\em cofoliation }on a stack. A cofoliation
is a change of the differentiable structure which amounts to giving a
full representable smooth epimorphism.  Cofoliations are uniquely
determined by their associated {\em Lie algebroids}.  

Cofoliations on stacks arise from {\em flat connections }on groupoids.
Connections on groupoids generalize connections on gerbes and bundles
in a natural way.  A flat connection on a groupoid is an integrable
distribution of the morphism space compatible with the groupoid
structure and complementary to both source and target fibres.  A
cofoliation of a stack determines the flat groupoid up to \'etale
equivalence.

We show how a cofoliation on a stack gives
rise to a refinement of the Hodge to De Rham spectral sequence, where
the $E_1$-term consists entirely of vector bundle valued cohomology
groups.  

Our  theory works for differentiable, holomorphic and
algebraic stacks. 
\end{abstract}
\tableofcontents

\section{Introduction}

This paper concerns the Hodge to de Rham spectral
sequence for stacks. By `stack' we mean `smooth Artin stack' in 
the differentiable, i.e., $C^\infty$-category, the holomorphic
category or the algebraic category over a field of characteristic
zero. 

\subsubsection{Manifolds}

Recall the Hodge to de Rham spectral sequence of a manifold $X$:
\begin{equation}\Label{mani}
E_1^{p,q}=H^q(X,\Omega^p)\Longrightarrow H_{DR}^{p+q}(X)\,.
\end{equation}
The $E_1$-term consists of the {\em Hodge cohomology groups}.  These
are the cohomology groups of the vector bundles $\Omega^p$, the
exterior powers of the cotangent bundle $\Omega_X$. The abutment is
the de Rham cohomology  of $X$.  

In the differentiable case, the $E_1$-term of the spectral sequence
(\ref{mani}) is just the de Rham  complex.  The same is true in the
holomorphic case if $X$ is Stein, or in the algebraic case if $X$ is
affine. 

The existence of the  spectral sequence (\ref{mani})
is almost trivial: by definition, the de Rham cohomology
$H^n_{DR}(X)=\hh^n\big(X,(\Omega\com,d)\big)$ is equal to the
hypercohomology of the de Rham complex, and the spectral sequence can
be obtained for example from the double complex one gets if one
calculates hypercohomology using \v Cech cochains.

Less trivial is the fact that $H^n_{DR}(X)$ is equal to the
`topological' cohomology $H^n(X)$ of $X$. In the differentiable and
holomorphic 
case this is the (differentiable or holomorphic) Poincar\'e lemma.  In
the algebraic case, it is a theorem of Grothendieck that the algebraic
de Rham cohomology 
$H^n_{DR}(X,\cc)=H^n(X^\an,\cc)$, the  `topological' cohomology of the
associated analytic manifold. 

Viewed in this light, the significance of (\ref{mani}) is that it
expresses the topological invariants $H^n(X)$ in terms of the
`coherent invariants $H^q(X,\Omega^p)$. By `coherent', we mean that
the $H^q(X,\Omega^p)$ are cohomology groups of coherent
$\O_X$-modules (in fact vector bundles). Of course, only the {\em
  objects }on the $E_1$-level are coherent, the differential comes
from non-coherent data, by which we simply mean that that the de Rham 
differential is not linear over functions (sections of $\O_X$).

\subsubsection{Stacks}

Any stack $X$ admits groupoid presentations $X_1\toto X_0$. To any
such groupoid presentation is associated a simplicial manifold
$X\lcom$.  In good 
cases, the de Rham cohomology of $X$ is defined as the total
cohomology of the double complex $\Gamma(X\lcom,\Omega\com)$, with \v
Cech and de Rham differentials (see, for example,~\cite{ictp}).

Alternatively, one may observe that there is a de Rham complex of big
sheaves $(\Omega_\bigsite\com,d)$ on the big site $X_\bigsite$ of the
stack $X$ and define 
$H^n_{DR}(X)$ as the hypercohomology of the big site of $X$ with
values in this big de Rham complex. (Note: the big site
  is defined as 
the category of all manifolds $U\to X$, with smooth structure map
$U\to X$.  No smoothness assumption on maps between various $U/X$ is
made. We use the \'etale topology on $X_\bigsite$.)

Either way, one obtains an $E_1$-spectral sequence
\begin{equation}\Label{false}
E_1^{p,q}=H^q(X_\bigsite,\Omega_\bigsite^p)\Longrightarrow
H_{DR}^{p+q}(X)\,, 
\end{equation}
abutting to the de Rham cohomology. It is also not difficult to reduce
the proof that $H_{DR}(X)$ is equal to the `topological' cohomology
of $X$ (i.e., the cohomology of the big topological site associated to
$X$) to the manifold case (see~[ibid.]). 

At first glance, the $E_1$-term of (\ref{false}) does not look very
coherent, as it involves big sheaves, which are by no means to be
considered coherent (much  less vector bundles).

One way to define the difference between big and small $\O_X$-modules
is as follows: If $\fF$ is a sheaf of $\O_X$-modules on the big site
of a stack $X$, then we get an induced sheaf of $\O_{X_n}$-modules
$\fF_n$ on every manifold $X_n$ partaking in the simplicial manifold
associated to a groupoid presentation of $X$. Moreover, we get
compatibility maps $\pi\upst \fF_m\to \fF_n$, for every structure
morphism $\pi:X_n\to X_m$. The sheaf $\fF$ is {\em coherent}, if every
$\fF_n$ is a coherent $\O_{X_n}$-module and all $\pi\upst \fF_m\to
\fF_n$ are isomorphisms.  For example, the big sheaf $\Omega_\bigsite$
induces $\Omega_{X_n}$ on $X_n$. Hence $\Omega_\bigsite$ is not
coherent, since only the first of the two coherence conditions is
satisfied.  We also define a vector bundle on a stack to be a
coherent $\O_X$-module $\fF$ all of whose components $\fF_n$ are
vector bundles.  Note that vector bundles have a well-defined rank,
something that a big sheaf such as $\Omega_\bigsite$ lacks.

The first indication that (\ref{false}) is not as bad as it looks, is
that for Deligne-Mumford stacks (in particular manifolds) we have 
$H^q(X_\bigsite,\Omega^p_\bigsite)=H^q(X,\Omega^p)$, and so
(\ref{false}) agrees with~(\ref{mani}).

\subsubsection{The cotangent complex}

The natural coherent analogue of the cotangent bundle of a manifold is
the cotangent complex $L_X$ of the stack $X$, which is an object of
the derived category of $\O_X$-modules. Because of our
smoothness assumption on all stacks, we can define it as the
homomorphism of big sheaves
$L_X=[\Omega_{\bigsite}\to \Omega_{\bigsite/X}]$. Here,
$\Omega_{\bigsite/X}$ is the big sheaf which induces on a smooth
$X$-manifold $U$ the sheaf $\Omega_{U/X}$. 

Note that $L_X$   is perfect, of
perfect amplitude contained in $[0,1]$.  This means  that $L_X$ is
locally quasi-isomorphic to the complex given by a homomorphism of
vector bundles.  In fact, if $X_1\toto X_0$  is a groupoid presenting
the stack $X$, then $L_X|_{X_0}=[\Omega_{X_0}\to\Omega_{X_0/X}]$.
(Proof: for any $U\to X_0$ we have 
$[\Omega_{X_0}|_U\to \Omega_{X_0/X}|_U]\stackrel{\scriptstyle\rm
  qis}{\longrightarrow}[\Omega_U\to\Omega_{U/X}]$.)
So once restricted to a presentation, the homomorphism
of big sheaves `collapses' to a homomorphism of vector bundles. (Note
that $\Omega_{X_0/X}$ can be identified with  the conormal bundle of
the identity section 
$X_0\to X_1$, so that $L_X|_{X_0}$ is the dual of
the Lie algebroid of the groupoid $X_1\toto X_0$.)

The  cotangent complex  is certainly `coherent' data, in the sense
that $L_X$ is an object of the derived category of $\O_X$-modules with
coherent cohomology.  

It is a remarkable fact, maybe first noticed by
Teleman (see \cite{tel}, for remarks on the equivariant case),  that
the natural 
homomorphism $L_X\to\Omega_\bigsite$, as well as its exterior powers,
induce isomorphisms on cohomology:
\begin{equation}\Label{tele}
\hh^q(X,\Lambda^p L_X)=H^q(X,\Omega^p_\bigsite)\,,
\end{equation}
for all $q,p$. Thus, we can rewrite the above spectral
sequence~(\ref{false}) as
\begin{equation}\Label{it}
E_1^{p,q}=\hh^q(X,\Lambda^p L_X)\Longrightarrow H^{p+q}_{DR}(X)\,.
\end{equation}
This looks much more like sequence (\ref{mani}) for manifolds.
Moreover, the objects on the $E_1$-level have now been shown to be
coherent data and the 
abutment is isomorphic to topological cohomology. We call (\ref{it})
the Hodge to de  Rham spectral sequence. 

\unsure{
There are two more requirements we would want to make on the spectral
sequence~(\ref{it}):
\begin{itemize}
\item (multiplicativity) The $E_1$-term is a bigraded algebra, where
the multiplication is given by combining  exterior product with cup
product.  We should require the differential on the $E_1$-level to be
a derivation.
\item (naturality) A morphism of stacks $X\to Y$ should induce a
  morphism of spectral sequences in the  other direction, respecting
  the differential graded algebra structures. 
\end{itemize}

Moreover, we might still include the following two examples:

(i) $B_MV$, for a vector bundle $V\to M$ over a manifold $M$.  Here it
seems that $E\to T_X$ is already the same thing as a connection on
$V$, at least if we require $M$ to give a pre-realization.

(ii) the case of a banal groupoid.  Here the spectral sequence induces
a flat connection on $\Upsilon$.  This cannot be encoded in $E\to T_X$
alone.}

\subsubsection{This paper}

So far the general story.  We now ask ourselves if we can still
`improve' upon~(\ref{it}), i.e., if we might be able to
construct a spectral sequence whose $E_1$-term actually consists of
cohomology groups of $X$ with values in vector bundles, as opposed to
just hypercohomology groups with values in perfect complexes.

We are motivated by the  case of a quotient stack $X=[X_0/G]$, for
which the cotangent complex is (quasi-isomorphic to) a homomorphism of
vector bundles, globally.  In fact, consider the  equivariant morphism
of equivariant 
vector bundles $\Omega_{X_0}\to\Gg^\vee$ on $X_0$. Here, $\Gg^\vee$ is the
trivial vector bundle associated to the dual
of the Lie algebra of $G$ with the adjoint representation,
$\Omega_{X_0}$ is the cotangent bundle of $X_0$ with the $G$-action
induced by the action of $G$ on $X_0$ and $\Omega_{X_0}\to\Gg^\vee$ is
obtained by differentiating the various orbit maps $G\to X_0$. 
Because $[\Omega_{X_0}\to\Gg^\vee]$ is $G$-equivariant,
it descends to a homomorphism of vector bundles $[\Omega\to \Upsilon]$
on the quotient stack $X$ and, in fact, the complex $[\Omega\to
\Upsilon]$ is quasi-isomorphic to 
the cotangent complex $L_{[X/G]}$. In other words, there is a distinguished
triangle
\begin{equation}\Label{pude}
\xymatrix@1{
L_X\rto & \Omega\rto & \Upsilon\rto & L_X[1]}
\end{equation}
in the derived category of $\O_X$-modules. (Note that pulling back
(\ref{pude}) via $\pi:X_0\to X$ gives the functoriality triangle
$\pi\upst L_X\to L_{X_0}\to L_{X_0/X}\to \pi\upst L_X[1]$.)

Since $\Lambda^p L_X=[\Omega^p\to\ldots\to S^p\Upsilon]$, one may ask
if there exists an $E_1$-spectral sequence whose $E_1$-term
consists of  the various $H^n(X,\Omega^{p-k}\otimes S^k\Upsilon)$ and which
abuts to $H_{DR}(X)$:
\begin{equation}\Label{wonder}
E_1^{m,n}=\bigoplus_{p+k=m}H^n(X,\Omega^{p-k}\otimes S^k\Upsilon)
\Longrightarrow H^{m+n}_{DR}(X)\,.
\end{equation}
In the equivariant case, such a spectral sequence has, indeed, been
constructed by Getzler \cite{equichern}. Rewriting~(\ref{wonder})
in equivariant language it reads:
\begin{equation}\Label{equilan}
E_1^{m,n}=\bigoplus_{p+k=m}H^n_G(X_0,\Omega_{X_0}^{p-k}\otimes
S^k\Gg^\vee) \Longrightarrow H^{m+n}_{DR,G}(X)\,.
\end{equation}
For example, if $G$ is compact (or reductive, in the algebraic
setting), then all higher $H^n_G$ vanish and we obtain a double
complex 
$$\bigoplus_{p,k}\Gamma(X_0,\Omega_{X_0}^{p-k}\otimes S^k\Gg^\vee)^G$$
computing the equivariant de Rham cohomology of $X$.  This is, of
course, nothing but the Cartan model. 

In this article, we examine to what extent a spectral
sequence~(\ref{wonder}) exists for general stacks. Note that
(\ref{wonder}) is also a generalization of (\ref{mani}).
But, in contrast with (\ref{it}), its $E_1$-term is computed in terms
of vector bundles.  

Moreover, one should think of (\ref{wonder}) as a refinement of
(\ref{it}). In fact (whenever we can construct (\ref{wonder}) at all),
there exists an
$E_1$-spectral sequence from (\ref{wonder}) to (\ref{it}).
For fixed $p$ it reduces on the $E_1$-level to the usual
$E_1$-spectral sequence 
\begin{equation}\Label{contemplate}
E_1^{k,n}=H^n(X,\Omega^{p-k}\otimes S^k\Upsilon)\Longrightarrow
\hh^{k+n}(X,\Lambda^p L_X)\,
\end{equation}
of hypercohomology.

\subsubsection{Flat connections on groupoids}

Let us describe what additional structure we need to have on our stack
$X$, for us to construct the spectral sequence~(\ref{wonder}). Of
course, we need a global resolution of the cotangent complex by vector
bundles $L_X=[\Omega\to\Upsilon]$. This global resolution
$[\Omega\to\Upsilon]$ needs to be endowed with an extra structure,
which we shall describe in the dual picture.

Giving a global resolution $L_X=[\Omega\to\Upsilon]$ of the cotangent
complex by vector bundles is equivalent to giving a surjective linear
map $a:E\to T_X$ from the vector bundle $E=\Omega^\vee$ to the tangent
stack $T_X$ of $X$. We define $N=\ker(E\to T_X)$, which can also be
constructed as the fibred product
\begin{equation}\Label{nfib}
\vcenter{\xymatrix{
N\rto\dto_\phi\ar@{}[dr]|\Box & X\dto^0\\
E\rto^a & T_X\,.}}
\end{equation}
It turns out that $N=\Upsilon^\vee$.

We define a {\em pre-realization }of $a:E\to T_X$ to be a presentation
$\pi:X_0\to X$ of the stack $X$, together with a vector bundle map
$p:T_{X_0}\to E$ covering $\pi$, such that $a\comp p$ is the canonical 
morphism $T_{X_0}\to T_X$ (i.e., there is a given 2-isomorphism, etc.),
and such that the rectangle marked with a box in
$$\xymatrix@C=1pc{
T_{X_0}\drto\ar[rrr]^p &&&
E\drto\rrto^a && T_X\dlto\\
& X_0\ar[rrr]^\pi \ar@{}[urr]|-\Box &&& X &}$$
is cartesian, i.e., a pullback diagram.  Such pre-realizations always
exist. 

It is easy to see that if $X_1\toto X_0$ denotes the groupoid induced
by the presentation $X_0\to X$ (which means that $X_1=X_0\times_X
X_0$), we get an induced subbundle $E_1\subset T_{X_1}$ defined by
$E_1=T_{X_0}\times_E T_{X_0}$.  (Use  the fact that
$T_{X_1}=T_{X_0}\times_{T_X} T_{X_0}$.)  Note also, that $E_1\toto T_{X_0}$
is a subgroupoid of $T_{X_1}\toto T_{X_0}$, and that $E_1$ is
complementary to both the source and target fibres of $X_1\toto X_0$.

We call such  a distribution $E_1\subset T_{X_1}$ a {\bf connection
}on the groupoid $X_1\toto X_0$. The key question is whether or not
the distribution $E_1$ is integrable: if we can find a pre-realization
of $a:E\to 
T_X$ for which the associated distribution $E_1\subset T_{X_1}$ is
integrable, we can construct the spectral sequence~(\ref{wonder}), at
least in good cases, for example in the algebraic case if the diagonal
of $X$ is affine. 

If the connection $E_1\subset T_{X_1}$ on a groupoid $X_1\toto X_0$ is
integrable as a distribution on $X_1$, we call the connection {\em
flat}. A groupoid endowed with a flat connection is called a {\em flat
  groupoid}. The terminology is justified by the relation to (flat)
connections on vector bundles and on gerbes. 

Flat connections on groupoids were independently discovered by 
Tang~\cite{XTang}. 

\subsubsection{Cofoliations}

In general, we have no means of comparing the spectral
sequences~(\ref{wonder}) coming from  different pre-realizations of
$a:E\to T_X$.  Thus we  have to put additional structure on $a:E\to
T_X$, to make the spectral sequence well-defined.\comment{there is no
  way to encode the connection on $\Upsilon$ purely in terms of $E\to
  T_X$}
We note that if the distribution $E_1\subset T_{X_1}$ is integrable,
it induces a structure on $a:E\to T_X$ comparable to that of {\em Lie
algebroid}.  If these induced  Lie algebroid structures on $a:E\to
T_X$ are the same, two different pre-realization give the same spectral
sequence~(\ref{wonder}). 

Thus  we are led to  study the Lie algebroid structures on  $a:E\to
T_X$, induced by flat groupoids pre-realizing $a:E\to T_X$.  For lack
of a better word\footnote{
we are open to alternative suggestions!}, 
we call these structures
{\bf cofoliations }on 
$X$. The spectral sequence~(\ref{wonder}) is thus an invariant of the
cofoliation on $X$.

A cofoliation is analogous to a foliation, in that we may think of a
cofoliation on $X$ as changing the differentiable (or holomorphic or
algebraic) structure on $X$. Except that the dimension goes up, instead
of down. For example, let $G$ be a Lie group and
$\widetilde{G}$ the same group but with discrete differentiable
structure. The quotient $[G/\widetilde{G}]$ is a cofoliation of the
point. It gives the point a differentiable structure of dimension
$\dim G$. 

In general, for a stack $X$, a cofoliation $\widetilde{X}\to X$ is a
change of differentiable structure making  $\widetilde{X}$ a
Deligne-Mumford stack and $\widetilde{X}\to X$ a full
smooth representable epimorphism.  (Analogously, a foliation of a
manifold is
a surjective immersion.)

Of course, this interpretation depends on accepting stacks with
horrible diagonal as geometric.  If one does not want to do this,
nothing is lost: this is just a way of thinking.  We do not actually use
stacks such as $[G/\widetilde{G}]$ in this article.

\subsubsection{The de Rham complex}

The spectral sequence~(\ref{wonder}) is constructed as follows: we
start with a flat groupoid $X_1\toto X_0$ representing and realizing
the stack $X$ with its cofoliation.  We obtain, among other things, a
global resolution $\phi:\Omega\to\Upsilon$ of $L_X$.  
We have
$$\Lambda^p(\Omega\to\Upsilon)=
\big(\bigoplus_{k=0}^{p}\Omega^{p-k}\otimes  
S^k\Upsilon,\phi\big)\,. $$
Then we take the
\v Cech complexes of
$\bigoplus_{p+k=m}\Omega^{p-k}\otimes S^k\Upsilon$, for all $m$,
associated to
the covering $X_0\to X$. Notation:
$$K^{p,k,n}=\Gamma(X_n,\Omega^{p-k}\otimes S^k\Upsilon)\,.$$
We have the \v Cech differential $\del$ on $K$, raising $n$ by 1. The
$\O_X$-linear map $\phi$ gives rise to a differential on $K$
increasing $k$ by 1.  We also construct, using flat connections on the
$\Upsilon|_{X_n}$, a differential $d$ on $K$, which increases $p$ by
1. Hence, the two differentials $\phi$ and $d$ both increase $m$ by
1.  We almost have, for every \v Cech degree $n$, a double complex
$(K^{\bullet,\bullet,n}, \phi+d)$, except for the fact that $\phi$ and
$d$ do not commute! To remedy the situation, we construct a homotopy
operator $\iota$ between the two compositions
$$\xymatrix@1{
(K^{p,k,\bullet},\del)\rrtwocell^{\phi\comp d}_{d\comp\phi}{\iota} &&
(K^{p+1,k+1,\bullet},\del)}\,.$$
Thus, after passing to \v Cech  cohomology, the two differentials
$\phi$ and $d$ commute, and we get, indeed, for every $n$, a double
complex
$$\big(H^n(X,\Omega^\bullet\otimes S^\bullet\Upsilon),\phi+d\big)\,.$$
The associated total complexes form the $E_1$-term of the spectral
sequence~(\ref{wonder}).  

In fact, $\iota$ is a fourth differential on
$K^{\bullet,\bullet,\bullet}$, and $(\phi+\del+d+\iota)$ is a
differential.  The cohomology of the total complex
$(K,\phi+\del+d+\iota)$ is the de Rham cohomology of the stack $X$
(thus the abutment of~(\ref{wonder}).)  The complex
$(K,\phi+\del+d+\iota)$ might be of independent interest.  Even though
it computes the same cohomology as the double complex
$\Gamma(X\lcom,\Omega\com)$ mentioned above, it is in some sense much
smaller:  it consists entirely of \v Cech cochains with values in
vector bundles (not big sheaves).

\subsubsection{Acknowledgements}

I would like to thank Ezra Getzler for useful discussions and for
referring me to his paper \cite{equichern}, which inspired this
research.

I would also like to thank Constantin Teleman for alerting me to
Formula~(\ref{tele}) and comments on an earlier version of this paper.

\unsure{
\subsubsection{To Do}

\begin{itemize}
\item
Prove isomorphism with Ezra's complex in equivariant case. 

\item
Proofs using coordinate systems in holomorphic and algebraic case.

\end{itemize}
}

\subsection{Notation}

\subsubsection{General setup}

We will be working in any of the following three categories: the
$C^\infty$-category (also called the differentiable category), the
holomorphic category, or the algebraic category, by which we mean the
category of finite type smooth schemes over a field of characteristic
zero.  (The characteristic zero assumption is important, because we
occasionally refer to the Lefschetz principle to reduce proofs to the
holomorphic case.) We use the word `manifold', to denote an object in
any of these categories. A `smooth map' is a morphism which is
pointwise surjective on tangent spaces, an `immersion' is a morphism
which is pointwise injective on tangent spaces. An \'etale map is a
smooth immersion.  (In the differential and holomorphic case, an
\'etale map is a local isomorphism.)

Let $\kk$ denote the real numbers if we are working in the
$C^\infty$-category, the complex numbers if we are in the holomorphic
category, or the ground field of characteristic zero if we are in the
algebraic category.

A manifold $X$ has a structure sheaf $\O_X$.  (In the
$C^\infty$-category this is the sheaf of $\rr$-valued differentiable
functions.)
The tangent bundle of a manifold $X$ is denoted by $T_X$, the
cotangent bundle and its exterior powers by $\Omega^p_X$. All vector
bundles are identified with their locally free sheaves of
$\O_X$-modules. 

\subsubsection{Groupoids and associated simplicial manifolds}

A {\bf groupoid} is a groupoid $X_1\toto X_0$ in our underlying
category whose source and target maps $s,t:X_1\to X_0$ are smooth.
The identity section of a groupoid is usually denoted by $\iota:X_0\to
X_1$. If
$s$ and $t$ are \'etale, we call the groupoid $X_1\toto X_0$ {\em
  \'etale}.  Every groupoid induces a `tangent groupoid' $T_{X_1}\toto
T_{X_0}$. 

If $\phi,\psi\in X_1$, such that $t(\phi)=s(\psi)$, we write the
composition as $\phi\ast\psi$. We also use  notation
$p_1(\phi,\psi)=\phi$ and  $p_2(\phi,\psi)=\psi$, as well as
$m(\phi,\psi)=\phi\ast\psi$. Right multiplication by $\psi$ is an
isomorphism 
$$R_\psi:t^{-1}\big(s(\psi)\big)\longrightarrow
t^{-1}\big(t(\psi)\big)\,,$$
with inverse $R_{\psi^{-1}}$. Similarly, left multiplication by
$\phi$ is an isomorphism
$$L_\phi:s^{-1}\big(t(\phi)\big)\longrightarrow
s^{-1}\big(s(\phi)\big)\,,$$
with inverse $L_{\phi^{-1}}$. 

To any groupoid $X_1\toto X_0$ is associated a simplicial manifold $X\lcom$.
We identify $X_n$ with the set of composable arrows
$$x_0 \stackrel{\psi_1}{\rightarrow}\ldots
\stackrel{\psi_{n}}{\rightarrow} x_{n}\,$$
in  $X_1\toto X_0$. 
For the structure maps of $X\lcom$ we use the following notation.
For $q=0,\ldots,n+1$ let
\[\hat{\pi}_q:{X}_{n+1}\longrightarrow {X}_n\]
be the projection leaving out the $q$-th object. In other words,
\begin{multline*}
\hat{\pi}_q:(x_0 \stackrel{\psi_1}{\rightarrow}\ldots
\stackrel{\psi_{n+1}}{\rightarrow} x_{n+1})
 \longmapsto\\
\longmapsto  (x_0 \stackrel{\psi_1}{\rightarrow}\ldots \rightarrow x_{q-1}
\stackrel{\psi_{q}\ast\psi_{q+1}}{\longrightarrow} x_{q+1}\rightarrow
\ldots \stackrel{\psi_{n+1}}{\rightarrow} 
x_{n+1})\,,
\end{multline*}
with obvious modifications for $q=0,n+1$.
For $q=0,\ldots,n-1$ we let
\[\iota_q:{X}_{n-1}\longrightarrow {X}_n\]
be the diagonal repeating the $q$-th object. In other words,
$$\iota_q:(x_0 \stackrel{\psi_1}{\rightarrow}\ldots
\stackrel{\psi_{n-1}}{\rightarrow} x_{n-1})\longmapsto
(x_0 \stackrel{\psi_1}{\rightarrow}\ldots \rightarrow x_{q}
\stackrel{\id}{\rightarrow} x_{q}\rightarrow \ldots
\stackrel{\psi_{n-1}}{\rightarrow} x_{n-1})\,.$$ 
We also denote by $\pi_q:X_n\to X_0$ the projection onto the $q$-th object:
$$\pi_q:(x_0 \stackrel{\psi_1}{\rightarrow}\ldots
\stackrel{\psi_{n}}{\rightarrow} x_{n})\longmapsto x_q\,,$$
$q=0,\ldots,n$. Finally, we have maps $\pi_{qr}:X_n\to X_1$, for
$0\leq q<r\leq n$, given by 
$$\pi_{qr}:(x_0 \stackrel{\psi_1}{\rightarrow}\ldots
\stackrel{\psi_{n}}{\rightarrow} x_{n})\longmapsto 
(x_q \stackrel{\psi_{q+1}\ast\ldots\ast\psi_r}{\longrightarrow}x_r)\,,$$

\subsubsection{Stacks}

In every one of our three categories we have a notion of stack.  We
will use the word {\em stack } for any stack over our base category
associated to a groupoid $X_1\toto X_0$ (which means that it is
isomorphic to the stack of torsors under $X_1\toto X_0$).  Since we are
only working with smooth objects $X_0$, $X_1$, all stacks are smooth. The stack
associated to the tangent groupoid $T_{X_1}\toto T_{X_0}$ is called
the {\em tangent stack}, denoted $T_X$. 

If the stack $X$ is given by the groupoid $X_1\toto X_0$, a {\em
vector bundle }over $X$ is any stack over $X$, induced by a groupoid
$E_1\toto E_0$ of vector bundles $E_i$ over $X_i$, where both diagrams
$$\xymatrix{E_1\rto\dto &  E_0\dto \\X_1\rto & X_0}$$
are cartesian, i.e., pullback diagrams.  Thus $T_X$ is {\em not }a
vector bundle over $X$, unless $X$ is Deligne-Mumford, which means we
can find an \'etale groupoid presenting $X$.

\subsubsection{The Lie algebroid of a groupoid}

Let us fix our conventions concerning the Lie algebroid of a groupoid
$X_1\toto X_0$. We denote by $N$ the normal bundle of $\iota:X_0\to
X_1$. The tangent bundle $T_{X_0}$, as well as $\iota\upst T_t$ and
$\iota\upst T_s$, are subbundles of $\iota\upst T_{X_1}$. Here $T_t$
and $T_s$ are the relative tangent bundles of the maps $s,t:X_1\to
X_0$. In fact, we have two direct sum decompositions
$$\iota\upst T_{X_1}=T_{X_0}\oplus \iota\upst T_t\qquad\text{and}\qquad
\iota\upst T_{X_1}=T_{X_0}\oplus \iota\upst T_s\,.$$
Any complement of $T_{X_0}$ is automatically isomorphic to $N$.
Thus there are two ways of thinking of $N$ as a subbundle of
$\iota\upst T_{X_1}$. There are also two projections from $\iota\upst
T_{X_1}$ onto $T_{X_0}$, one with kernel $\iota\upst T_t$ (denoted
$\eta_t$) and one with
kernel $\iota\upst T_s$ (denoted $\eta_s$.). Composing, we get four
maps $N\to T_{X_0}$, 
two of which vanish.  Unfortunately, the other two are often not equal, but,
rather, add up to zero.  Thus we have to make up our mind, which of
these two maps we declare to be the {\em canonical }map $\phi:N\to
T_{X_0}$. (This map is also called the `anchor map'.) We make our
choice as follows: \comment{I'm not happy with this part yet}
$$\text{\em We identify $N$ with $\iota\upst T_t$ using the
composition $\iota\upst T_t\to\iota\upst T_{X_1}\to N$,}$$
and then
apply $\eta_s$. This makes $\phi$ essentially equal to $\eta_s$. We
will also need to declare a 
`canonical' isomorphism $\iota\upst T_s\iso N$.  We choose the {\em
  negative }of the composition $\iota\upst T_s\to T_{X_1}\to N$.  Thus
we have declared two canonical injections $N\to \iota\upst
T_{X_1}$. We denote the inclusion, with image $\iota\upst T_t$, by
$\rho_s:N\to\iota\upst T_{X_1}$, and the other injection by $\rho_t$. The
way we made our choices, we have
$\eta_t\comp\rho_t=\eta_s\comp\rho_s=\phi$ (and
$\eta_t\comp\rho_s=\eta_s\comp\rho_t =0$). 

Now that we have chosen $\phi:N\to T_{X_0}$, we can think of $N$ as the
relative tangent bundle of the map $X_0\to X$, where $X$ is the stack
associated to $X_1\toto X_0$. To justify this, let us define
identifications $s\upst N\iso T_t$ and $t\upst N\iso T_s$. Fix a
point $\psi\in X_1$. We have $s\upst
N|_\psi=N|_{s(\psi)}=T_t|_{s(\psi)}$. Then we compose with the
derivative of right multiplication by $\psi$:
$$DR_\psi|_{s(\psi)}:T_t|_{s(\psi)}\longrightarrow
T_t|_\psi\,$$
to obtain the identification $s\upst N=T_t$. Similarly, we have
$t\upst N|_\psi=N|_{t(\psi)}\cong T_s|_{t(\psi)}$. Now composing with the
derivative of left multiplication by $\psi$:
$$DL_\psi|_{t(\psi)}:T_s|_{t(\psi)}\longrightarrow
T_s|_\psi\,,$$
we obtain the identification $t\upst N=T_s$. 
It is important to notice that the two compositions
$$p_1\upst t\upst N\longiso p_1\upst T_s\longiso T_m\qquad\text{and}
\qquad
p_2\upst s\upst N\longiso p_2\upst T_t\longiso T_m$$
are equal.\comment{please check this. It is very important.}
This implies that we have canonical identifications $\pi_q\upst N\cong
T_{\hat\pi_q}$, for all $q=0,\ldots,n$ on every $X_n$.

\section{Flat Connections on Groupoids}

\subsection{Definition}

Recall that an integrable distribution on a manifold $X$ is a subbundle
$E\subset T_X$ of the tangent bundle such that $E^\perp\subset\Omega_X$
generates an ideal in $\bigoplus_p\Omega^p_{X}$ which is preserved by the
exterior derivative. Equivalently, the sheaf $E$ is closed under the
Lie bracket inside $T_X$.

\begin{defn}\Label{parst}
A {\bf connection }on the groupoid $X_1\toto X_0$ is a subbundle
$E\subset T_{X_1}$, such that $E\toto T_{X_0}$ is a subgroupoid 
of $T_{X_1}\toto T_{X_0}$ and both 
diagrams
$$\xymatrix{
E\ar@<-.5ex>[r]\ar@<.5ex>[r]\dto & T_{X_0}\dto\\
X_1\ar@<-.5ex>[r]\ar@<.5ex>[r] & X_0}$$
are cartesian.

The connection $E$ is {\bf flat }or {\bf integrable}, if $E\subset
T_{X_1}$ is an integrable distribution.
\end{defn}

\begin{defn}
A groupoid endowed with an integrable connection, will be called a
{\bf flat groupoid}. 
\end{defn}

\begin{rmk}
Flat connections on groupoids are called {\em  \'etalifications
}in~\cite{XTang}. 
\end{rmk}

\begin{numrmk}
In the differential and holomorphic categories, the integrable
distribution $E$ admits integral submanifolds.\comment{reference?}
Denote by $\tilde{X}_1$ the union of the leaves of this
foliation. Note that $\tilde{X}_1$ has the same set of points as
$X_1$, but a different differentiable (holomorphic) and topological
structure. The canonical map $\tilde{X}_1\to X_1$ is a bijective
immersion.  The conditions on the flat connection $E$ are equivalent
to saying that $\tilde{X}_1\toto X_0$ is an \'etale groupoid and that
$\tilde{X}_1\to X_1$ is a morphism of groupoids.
\end{numrmk}

We shall denote the groupoid $E\toto T_{X_0}$ also by $E_1\toto
E_0$. As usual, we get a simplicial manifold $E\lcom$, which is, in
fact, a simplicial submanifold of $T_{X\lcom}$.  Note  that
$E\lcom$ is a vector bundle over $X\lcom$ (which $T_{X\lcom}$ is {\em
not\/}).  If $E$ is a flat connection, the subbundle  $E_n\subset
T_{X_n}$ is an integrable distribution on $X_n$, for all $n\geq0$.

Let $N_0=N_{X_0/X_1}$ be the normal bundle to the identity of
$X_1\toto X_0$. We always identify $N_0$ with $T_{X_0/X}$, the
relative tangent bundle of $X_0$ over the stack $X$ defined by the
groupoid $X\lcom$. The relative tangent bundle of $\hat{\pi}_q:X_n\to
X_{n-1}$, is canonically identified with
$\pi_q\upst N_0$. Thus, for all $q=0,\ldots,n$, we have a canonical subbundle
$\pi_q\upst N_0\subset T_{X_n}$ (which is an integrable distribution).

If $E$ is a connection, then for every $q=0,\ldots,n$, we have a
direct sum decomposition 
\begin{equation}\Label{dirsu}
T_{X_n}=E_n\oplus\bigoplus_{r\not=q} \pi_r\upst N_0\,.
\end{equation}
In particular, $T_{X_1}=E\oplus s\upst N_0=E\oplus t\upst N_0$. These
splittings give rise to projections $\omega:T_{X_1}\to s\upst N_0$ and
$\tilde{\omega}:T_{X_1}\to t\upst N_0$ and an isomorphism $t\upst N_0\to s\upst
N_0$ compatible with the projections from $T_{X_1}$. Let us call
$\omega\in\Gamma(X_1,\Omega_{X_1}\otimes s\upst N_0)$ the {\bf differential
form }of the  connection $E$. Of course, $E$ can be recovered from
$\omega$ as its kernel. Note that the isomorphism $t\upst N_0\to s\upst N_0$ is
nothing but the restriction of $\omega$. 

\begin{lem}\comment{need to characterize integrability in terms of $\omega$.}
The differential form $\omega:T_{X_1}\to s\upst N$ of a connection
$E$ on $X_1\toto X_0$ satisfies:

(i) $\omega\resto s\upst N=\id_{s\upst N}$,

(ii) $\iota\upst\omega:T_{X_0}\to N$ satisfies $\iota\upst\omega=0$,

(iii) as homomorphisms $T_{X_2}\to\pi_0\upst N$, we have
$\mu\upst\omega=p_1\upst\omega+\chi_{01}(p_1\upst\omega)$. 
\end{lem}

\unsure{
Conversely, if $\omega:T_{X_1}\to s\upst N$ satisfies these three
properties, then $E=\ker\omega$ is a connection on $X_1\toto X_0$
with associated form $\omega$. Except, this does not really make
sense, because $\chi_{01}$ is defined in terms of $\omega$.
}

Using the decomposition (\ref{dirsu}) for $q=q_1$ and $q=q_2$, where
$q_1<q_2$, gives us an isomorphism $\pi_{q_2}\upst N_0\to \pi_{q_1}\upst N_0$,
whose {\em negative }we shall denote by $\chi_{q_1,q_2}$. Thus we have an
anti-commutative diagram 
$$\xymatrix{ {T_{X_n}}\rto\drto & {\pi_{q_2}\upst
N_0}\dto^{\chi_{q_1,q_2}} \\ & {\phantom{\,.}\pi_{q_1}\upst N_0\,.}}$$ Note
that we may identify $\chi_{q_1,q_2}$ with $-\pi_{q_1,q_2}\upst\omega$.

The canonical projection $T_{X_n}\to\pi_q\upst T_{X_0}$ induces an
isomorphism $E_n\to\pi_q\upst T_{X_0}$, and hence, by composition,
isomorphisms $\psi_{q_1,q_2}:\pi_{q_2}\upst T_{X_0}\to\pi_{q_1}\upst
T_{X_0}$, for any $q_1$, $q_2$. Denoting the canonical homomorphism
$N_0\to T_{X_0}$ by $\phi$, the diagrams
$$\xymatrix{
{\pi_{q_2}\upst N_0}\dto_{\chi_{q_1,q_2}}\rrto^{\pi_{q_2}\upst\phi} &&
{\pi_{q_2}\upst T_{X_0}}\dto^{\psi_{q_1,q_2}}\\
{\pi_{q_1}\upst N_0}\rrto^{\pi_{q_1}\upst\phi} &&
{\pi_{q_1}\upst T_{X_0}}}$$
commute. In other words, we have defined {\em descent data }for the
homomorphism of vector bundles $\phi:N_0\to T_{X_0}$. Let us denote the
induced homomorphism of vector bundles on $X$, the stack associated to
$X_1\toto X_0$, by $\phi:N\to E$. We call $N\to E$ the {\bf
  representative of the tangent complex }given by our
connection. Let us denote the pullback of $\phi:N\to E$ to $X_n$ by 
$\phi:N_n\to E_n$. 

\begin{rmk}
Another way to construct $\phi:N\to E$ is as follows:  the groupoid
morphism $E\lcom\to T_{X\lcom}$ gives rise to a linear epimorphism of
$X$-stacks $E\to T_X$. Let $N$ be the kernel, constructed as in
Diagram~(\ref{nfib}).  Then we have 
$$T_X=[E/N]\,,$$
the quotient (over $X$) of the vector bundle $E$ by the action of the
vector bundle $N$ by addition via $\phi$.

Duality between the tangent stack and the
cotangent complex implies that we have constructed a distinguished
triangle of complexes of $\O_X$-modules
$$\xymatrix@1{L_X\rto & E^\vee \rto & N^\vee \rto & L_X[1]}\,.$$
We will use notation $\Omega=E^\vee$ and $\Upsilon=N^\vee$.
\end{rmk}

\unsure{
Note that we have, in fact, a
distinguished triangle of complexes of $\O_X$-modules
$$N\stackrel{\phi}{\longrightarrow} E \longrightarrow T_X \longrightarrow
N[1]\,$$
\comment{I'm confused about the homo $T_X\to N[1]$  Also, I don't
  think $T_X$ is an $\O_X$-module!!!}
where $T_X$ is the tangent complex of $X$.
}

\begin{rmk}
A connection on $X\lcom$ is the same thing as descent data for the
distinguished triangle
$$\xymatrix@1{\pi\upst L_X\rto & L_{X_0}\rto & L_{X_0/X}\rto & 
\pi\upst L_X[1]}\,$$
inducing the trivial descent data on $\pi\upst L_{X}$.
\end{rmk}

\subsubsection{More notation}

Now, for every $q$, we have an identification $N_n=\pi_q\upst
N_0$, and hence $n+1$ canonical ways of thinking of $N_n$ as a
subbundle of $T_{X_n}$. We denote these various embeddings by
$\rho_q:N_n\hookrightarrow T_{X_n}$. 

For every choice of $0\leq q < r\leq n$, we also get a canonical way
of making $N_n$ into a quotient of $T_{X_n}$. In fact, take
$$T_{X_n}\to \pi_{qr}\upst
T_{X_1}\stackrel{\pi_{qr}\upst\omega}{\longrightarrow} \pi_q\upst N
=N_n\,.$$
Let us denote this quotient map by $\omega_{qr}:T_{X_n}\to N_n$. We
have that $\omega_{qr}(E_n)=\omega_{qr}(\pi_j\upst N)=0$, for all
$j\not=q,r$. Moreover, $\omega_{qr}(\pi_q\upst N)=\id$, and
$\omega_{qr}(\pi_r\upst N)=-\id$. For $q>r$, let us define
$\omega_{qr}=-\omega_{rq}$, and for $q=r$, set $\omega_{qr}=0$. 

We also have, via the various identifications $E_n=\pi_q\upst T_{X_0}$,
for every $q$ a way of considering $E_n$ as a quotient of $T_{X_n}$. Let
us denote the corresponding quotient map by $\eta_q:T_{X_n}\to E_n$. 

Note that for every $q,r=0,\ldots,n$ the diagram
\begin{equation}\Label{diaqr}
\vcenter{\xymatrix{
N_n\drto_{\delta_{qr}\phi}\rto^{\rho_q} & T_{X_n} \dto^{\eta_r} \\
& E_n}}
\end{equation}
commutes, where $\delta_{qr}$ is the delta of Kronecker.

For future reference, let us also remark  the commutativity of the
diagram
\begin{equation}\Label{sumdiagram}
\vcenter{\xymatrix{
N_n\rto^\phi\drto_{\sum_{q=0}^n\rho_q} & E_n\dto\\ & T_{X_n} }}
\end{equation}

For the convenience of the reader, let us summarize the various maps
we constructed:

Without a connection, we have canonical maps 
$$\pi_q\upst N\stackrel{\rho_q}{\longrightarrow}T_{X_n}
\stackrel{\eta_r}{\longrightarrow} \pi_r\upst T_{X_0}\,,$$
but we use the connection to identify all $\pi_q\upst N$ with
$N_n$ and all $\pi_r\upst T_{X_0}$ with $E_n$. 

by adding the connection, we get maps
$$E_n\stackrel{\delta}{\longrightarrow}T_{X_n}
  \stackrel{\omega_{qr}}{\longrightarrow}    N_n\,.$$

\subsubsection{A few functorial properties}

\begin{defn}
Let $f\lcom:X\lcom\to Y\lcom$ be a morphism of groupoids. Suppose
$X\lcom$ and $Y\lcom$ are endowed with connections $E\subset
T_{X_1}$ and $F\subset T_{Y_1}$. The morphism $f$ is called  {\bf
  horizontal}, if the derivative $D{f_1}:T_{X_1}\to 
T_{Y_1}$ maps $E$ into $F$. 
\end{defn}

\begin{rmk}
Given two horizontal morphisms $f,g:X\lcom\to Y\lcom$, there is an
obvious notion of {\em horizontal natural transformation }from $f$ to
$g$. 
\end{rmk}

A morphism of groupoids $f\lcom:X\lcom\to Y\lcom$ is called {\em
\'etale}, if both $f_0:X_0\to Y_0$ and $f_1:X_1\to Y_1$ are
\'etale.  Suppose $X\lcom$ and $Y\lcom$ are endowed with connections
$E$ and $F$, respectively.  If $f$ is \'etale and 
horizontal, then the diagram
\begin{equation}\Label{ethor}
\vcenter{\xymatrix{E\rto\dto &F\dto\\
T_{X_1}\rto & T_{Y_1}}}
\end{equation}
is cartesian.

Conversely, given an \'etale morphism $f:X\lcom\to Y\lcom$ of
groupoids, any connection $F$ on $Y\lcom$ induces a unique connection
$E$ on  $X\lcom$ making $f$ horizontal: define $E$ as the fibred
product~(\ref{ethor}). 
We write $E=f\upst F$. If $F$ is integrable, then so is $E$.

A morphism of groupoids $f\lcom:X\lcom\to Y\lcom$ is called {\em
  cartesian}, if both diagrams
$$\xymatrix{
X_1\ar@<-.5ex>[d]\ar@<.5ex>[d]\rto^{f_1} &
Y_1\ar@<-.5ex>[d]\ar@<.5ex>[d]\\ 
X_0\rto^{f_0} & Y_0}$$
are cartesian. If $F$ is a connection on $Y\lcom$ and $f$ is
cartesian, we get an induced connection $E$ on $X\lcom$ by setting
$E=F\times_{T_{Y_1}} T_{X_1}$. If $F$ is integrable, so is $E$.

\subsection{The derived connection}

A connection $E\subset T_{X_1}$ on a groupoid $X_1\toto X_0$
induces a connection $\nabla$ on the vector bundle $N_0$, as follows. We
define, for 
$v\in\Gamma(X_0,T_{X_0})$ and $\nu\in\Gamma(X_0,N_0)$ the covariant
derivative of $\nu$ with respect to $v$ by
$$\nabla_v(\nu)=\iota\upst\omega[\rho_s(s\upst\nu), s\upst(v)]\,,$$
where $[\cdot\,,\cdot\,]$ denotes the Lie bracket of vector fields on
$X_1$. We identify $T_{X_0}$ with $E_0$, so that $s\upst(v)$, which is
a section of $E_1$, is canonically a section of $T_{X_1}$. The map
$\rho_s$ is the inclusion map $s\upst N_0\hookrightarrow T_{X_1}$. 
Note that the `mirror' formula 
$$\widetilde{\nabla}_v(\nu)=\iota\upst\tilde{\omega}[\rho_t(t\upst\nu),
t\upst(v)]\,$$ 
gives rise the the same connection, $\widetilde{\nabla}=\nabla$. We call
$\nabla$ the {\bf derived  connection }associated to the connection
$E$. If $E$ is flat, $\nabla$ is integrable.

\begin{numrmk}\comment{this is not so deep}
Note the following: given a  smooth map of manifolds $f:X\to Y$, a
distribution $E\hookrightarrow T_X$, such that  $E\iso
f\upst T_Y$, and a section $\iota:Y\to X$ whose image is a leaf (i.e.,
parallel to $E$), we get an
induced connection on $\iota\upst T_{X/Y}=N_{Y/X}$, where we
think of $Y$ is a submanifold of $X$ via $\iota$. 

Applying this principle to $s:X_1\to X_0$ with the identity section,
we get a connection on $N_0$.  We get the same connection
applying this principle to $t:X_1\to X_0$. This gives us a more
geometric way of defining the derived connection $\nabla$. 
\end{numrmk}
 
\begin{numrmk}
For every groupoid $X_1\toto X_0$ we get an associated {\em derived
}groupoid $N\toto X_0$, which is simply the normal bundle to the
identity section considered as a family of Lie groups.  This process
of deriving a groupoid is functorial and commutes with passing to
tangent groupoids.  Thus, if we derive the diagram 
$$\xymatrix@1{E\lcom\rto & T_{X\lcom}\rto & X\lcom}$$
of groupoids, we obtain the sequence
$$\xymatrix{
N_{E_0/E_1}\dto\rto &T_{N}\rto\dto & N\dto\\
E_0\rto & T_{X_0}\rto & X_0}$$
of vector bundles.  In other words, we have a {\em derived }connection
on the derived groupoid $N\to X_0$. A groupoid connection on a vector
bundle considered as a groupoid is the same thing as a vector bundle
connection in the traditional sense (see below).  Thus, this process
of deriving a 
groupoid induces a connection on the vector bundle $N\to X_0$.  This
connection is the derived connection defined above.
\end{numrmk}

Pulling back to $X_1$, we get connections $s\upst\nabla$ on $s\upst N_0$ and
$t\upst\nabla$ on $t\upst N_0$. Note that, in general, the isomorphism
$\chi:t\upst N_0\to s\upst N_0$ is not horizontal. Thus, the
connection $\nabla$ on $N_0$ does not descend to a connection on $N$ over
the stack $X$.

Using $\chi$ to identify $s\upst N_0$ with $t\upst N_0$, we get two
different flat connections on $N_1$, denoted $s\upst \nabla$ and
$t\upst \nabla$.  Their difference $\Psi_1=s\upst\nabla-t\upst\nabla$ is
a vector bundle homomorphism $\Psi_1:T_{X_1}\otimes N_1\to N_1$. 

\begin{lem}\Label{SameReason}
For the
same reason that $\widetilde\nabla=\nabla$, we have that $\Psi_1$
vanishes on $E\otimes N_1$:  the two connections $s\upst\nabla$ and
$t\upst\nabla$ agree on $E\otimes N_1$.\qed
\end{lem}

\begin{rmk}
The quotient $T_{X_1}/E$ being canonically isomorphic to $N_1$, we see
that $\Psi_1$ induces a vector bundle homomorphism $\Psi_1:N_1\otimes
N_1\to N_1$. Restricting back (via $\iota$) to $X_0$, we get a vector
bundle homomorphism $\Psi_0:N_0\otimes N_0\to N_0$.  This pairing on
$N_0$ agrees with the Lie algebroid pairing on $N_0$ if we restrict to
covariantly constant sections of $N_0$. 
\end{rmk}

\begin{numrmk}\Label{afterSR}
Passing to $X_n$, we get $n+1$ different connections on $N_n$, each
one given by one of the identifications $N_n=\pi_q\upst N_0$. In
other words, all $n+1$ subbundles $\eta_q:N_n\hookrightarrow T_{X_n}$
are endowed in a canonical way with an integrable
connection.  Let us denote the connection $\pi_q\upst\nabla$ on $N_n$
by $\nabla^q$. Then, for $0\leq q<r\leq n$, the difference
$\nabla^q-\nabla^r$ is equal to the composition
$$N\stackrel{\pi_{qr}\upst\Psi_1}\longrightarrow \pi_{qr}\upst
N^\vee\otimes
N\stackrel{\omega_{qr}}{\longrightarrow} \Omega_{X_n}\otimes N\,,$$
or, more succinctly, $\nabla^q-\nabla^r=\omega_{qr}(\Psi)$, which now
holds for all $q,r$. The most important consequence of these
considerations for us is that $\nabla^q-\nabla^r$ takes values in
$E_n^\perp\otimes N$, for all $q,r$.
\end{numrmk} 

\begin{rmk}\comment{Improve the treatment of $\Psi$}
Note that $s\upst \Psi_0= t\upst\Psi_0=\Psi_1$, so that $\Psi_0$
descends to $X$: there is a homomorphism of vector bundles $N\otimes
N\to N$ on $X$, making $N$ a bundle of Lie algebras on the stack
$X$.\comment{I only convinced myself in the equivariant case}
See also Remark~\ref{Liebundle}.
\end{rmk}

\subsection{Examples}

\begin{numex}
An \'etale groupoid  has a unique connection.  It is integrable.
\end{numex}

\begin{numex}
Let $G$ be a Lie group, considered as a groupoid $G\toto\ast$. There
is a unique connection on this groupoid, namely
$0\hookrightarrow T_G$.  It is integrable.  Of course, $N_0=\Gg$, the
Lie algebra of 
$G$. The derived connection is the unique (trivial) connection on the
vector space $\Gg$.  The descent datum 
$\chi:\Gg_X\to\Gg_X$ is the adjoint representation.  It is not locally
constant, unless $G$ is abelian or discrete (or more generally, has an
abelian connected component). This provides examples where the
derived connection does not descend to the stack.
\end{numex}

\begin{numex}
If  $X_1\toto X_0$ is a transformation groupoid $G\times X\toto
X$, letting $E=T_{G\times X/G}$ defines the {\em canonical }flat connection.
\end{numex}

\begin{numex}
For a manifold $X$, a connection on the groupoid $X\times X\toto X$
is the same thing as a trivialization of the tangent bundle $T_X$ of
$X$. (N.B. By a trivialization we mean trivialization up to choice of
basis: more precisely, descent data of $T_X$ to the point.) 
 
The connection is flat if and only if the Lie bracket is constant.
This means that there exists a Lie algebra $\Gg$ and an identification
$T_X=\Gg_X$ such that for {\em constant }sections of $T_X$ the Lie
bracket as vector fields coincides with the Lie bracket as elements of
$\Gg$. 
\end{numex}

\begin{numex}\Label{banalex}
More generally, if $f:X\to Y$ is a surjective submersion, and $X_1\toto X_0$
is the associate {\em banal }groupoid given by $X_0=X$ and
$X_1=X\times_YX$, a connection on $X_1\toto X_0$ is the same thing
as a triple 
$(E,\epsilon,\psi)$, where $\epsilon:E\to T_Y$ is an epimorphism of vector
bundles on $Y$ and $\psi:T_X\to E$ is a homomorphism of vector bundles
covering $f$, identifying $T_X$ with $f\upst E$ and such that
$\epsilon\comp\psi=Df$.
Less formally, we may say that a connection is descent data
for $T_X$, compatible with the map to $f\upst T_Y$.  

We will see later (Lemma~\ref{banalem}) what flatness means in this
context. 
\end{numex}

\begin{numex}
If $F\to X$ is a vector bundle over a manifold $X$ and we consider the
groupoid $X_1\toto X_0$ given by $X_0=X$ and $X_1=F$, with the groupoid
structure given by vector addition, a groupoid connection on $X_1\toto X_0$
is the same thing as a vector bundle connection on $F$. Most directly,
this can be 
seen by using the characterization of vector bundle connections in terms of
horizontal distributions on the total space of $F$. Also, we may
remark that $F=N_{X_0/X_1}$ and the
derived  connection determines the groupoid connection.  Integrability
of the groupoid connection is equivalent to integrability of the vector
bundle connection.
\end{numex}

\begin{numex}
More generally, if $X_1\toto X_0$ is a family of groups $G\to X$, a
connection identifies infinitesimally close fibres of $G\to X$ with
each other as
group schemes. For example, for a family of elliptic curves (in the
holomorphic or algebraic context), this means that the $j$-invariant of the
family is locally constant.  In this example, the absence of
infinitesimal automorphisms of elliptic curves implies that every
connection is flat.
\end{numex}

\subsubsection{Relation to connections on gerbes}

Let $X_1\toto X_0$ be an \'etale groupoid and $R_1\toto X_0$ a
groupoid with the same base space endowed with a morphism of groupoids
$R\lcom\to X\lcom$. Assume that $R_1\to X_1$ is a surjective
submersion.
(This data induces a morphism of stacks $\GG\to
X$, where $\GG$ is a {\em gerbe }over the Deligne-Mumford stack $X$.)
The kernel of $R_1\to X_1$ is a family of groups $\gG\to X_0$ over
$X_0$. Any connection on $R\lcom$ induces one on $\gG\toto
X_0$. (Note that $X_1\toto X_0$ being \'etale, the diagonal $X_0\to
X_1$ is an open immersion, hence the same is true for $\gG\subset
X_1$.  The restriction of the connection from $X_1$ to $\gG$
is a special case pulling back via \'etale morphisms, discussed
above.)

Suppose that we have given an identification of $\gG\to X_0$ with a
constant family of groups $G\times X_0\to X_0$, where $G$ is a Lie
group, and that $G\times X_0$ is central in $R_1$.  In this case, the
gerbe $\GG\to X$ is a $G$-gerbe and $R\lcom\to X\lcom$ is a 
$G$-central extension of groupoids. Note that $R_1\to X_1$ is now a
principal $G$-bundle, the $G$-action being induced from the groupoid
multiplication.  Moreover, the product family $G\times X_0\toto X_0$ has a
canonical flat connection on it.

Assuming given a connection $E$ on $R\lcom$, which restricts to the
canonical one one $G\times X_0\subset R_1$,
the compatibility of $E$ with the groupoid multiplication implies the
invariance of $E$ under the $G$-structure.  Thus, $E$ is a bundle
connection on the $G$-bundle $R_1\to X_1$. Let $\theta$ be the
connection form (which may be identified with the form $\omega$ of the
groupoid connection). One checks that the \v 
Cech coboundary of $\theta$ vanishes.  Thus we have defined a
{\em connective structure }on the groupoid central extension
$R\lcom\to X\lcom$, and hence on the $G$-gerbe $\GG$. (See 
\cite{BeXu}, where the case of $G=S^1$ is treated. See also the
classic reference for connections  on gerbes~\cite{brylinski}.)

If the groupoid connection $E$ is flat, then the connective structure
is  flat (more precisely, we may take $0$ as curving and then the
curvature vanishes). Thus a flat connection on the groupoid $R\lcom$,
extending the canonical one on $G\times X_0\subset R_1$, 
induces a flat connection on the $G$-gerbe $\GG$.

Conversely, given a $G$-gerbe $\GG\to X$ over a Deligne-Mumford stack,
we can find a groupoid central extension $R\lcom\to X\lcom$ 
with kernel $G\times X_0$, inducing $\GG\to X$.  We do this by
choosing a local trivialization of $\GG$ over an \'etale cover $X_0\to
X$.  A connective structure on the $G$-gerbe $\GG$ then defines a
connection on the groupoid $R\lcom$, which restricts to the canonical
one on the kernel, if the local trivialization $X_0\to \GG$ giving
rise to the central extension is compatible with the connective
structure.

If $\GG$ is endowed with a flat connection, then, at least in the
differentiable and holomorphic context, we can locally trivialize
$\GG$ as a flat gerbe. Then the induced connection on the induced
groupoid central extension $R\lcom\to X\lcom$ will be flat.  
(If we use
a  local trivialization which respects only the connective structure
on the flat gerbe to write $\GG$ is a $G$-central
extension, we get a connection on $R\lcom$ which is not flat.  But we
will  get a covariantly closed 2-form $B$ with values in $N$ on $X_0$,
whose \v Cech coboundary is equal to the curvature of the bundle
connection on $R_1$. The form $B$ represents the incompatibility of
the local trivialization with the {\em
curving\/}).

\begin{numrmk}\Label{conceivable}
It is conceivable, that, in general, there exists a notion of a {\em
trivialization of the curvature }for a connection on a groupoid.  This
would be a 2-form $B\in 
\Omega_{X_0}^2\otimes N$, such that $t\upst B-s\upst B$ is the
curvature of the distribution $E_1\subset T_{X_1}$. (Recall that the
curvature is a tensor $\Lambda^2 E_1\to T_{X_1}/E_1=N_1$.) Then the
existence of a covariantly closed $B$ might  serve as a
generalization or substitution of flatness of groupoid connections.
\end{numrmk}

\unsure{This concept might be useful especially in the algebraic context,
where cofoliations are probably less likely to admit realizations.
Working with pre-realizations, $B$ might do the job.

Can one still construct the spectral sequence?  Now $d$ will be a
differential at most up to homotopy.}

\section{Cofoliations}

Let us make a few remarks about Lie algebroids.
Note that if $a$ is
surjective, its kernel is a vector bundle, and in fact a bundle of Lie
algebras (i.e., the bracket is linear over functions). 

\begin{defn}
Let $f:V\to U$ be a morphism of manifolds.  Suppose $F\to T_V$ and $E\to
T_U$ are Lie algebroids. We call a morphism $\phi:F\to E$ a {\bf
  homomorphism }of Lie algebroids covering $f$, if 

(i) $\phi$ is a
homomorphism of vector bundles covering $f$, 

(ii) $\phi$ commutes with the
anchor maps, 

(iii) the brackets are
compatible, which means that if $s_0\in\Gamma(U,E)$ and 
$t_0\in\Gamma(V,F)$ form a compatible pair of sections, i.e.,
$f\upst(s_0)=\phi\lst(t_0)$,  and $s_1$, $t_1$ form another such pair,
then $[s_0,s_1]_E$ and $[t_0,t_1]_F$ are compatible in the same way.
\end{defn}

\subsection{Cofoliations on Manifolds}

Over a manifold $X$, an {\em integrable distribution }is essentially
the same thing as a Lie algebroid $a:E\to T_X$, whose anchor map $a$
is a subbundle.  We call an integrable distribution a {\em foliation},
if there exist submanifolds integrating $E\to
T_X$. Of course, in the differentiable or holomorphic context, every
integrable distribution is a foliation, by the Frobenius integrability
theorem.  But  in the algebraic category, there is a difference.

If,
instead, the anchor of a Lie algebroid is surjective (as a
homomorphism of vector bundles) we suggest calling the induced
structure on $X$ a {\em cofoliation}, at least if the  analogues of
integral submanifolds exist.

\begin{defn}\Label{realize}
Let $a:E\to T_X$ be a Lie algebroid on the manifold $X$. A {\bf
local  realization }of $E$ is given by the following data:

(i) a manifold $U$,

(ii) a smooth map $\pi:U\to X$,

(iii) a vector bundle homomorphism $p:T_U\to E$ covering $\pi$.

\noindent This data is subject to the constraints:

(iv) the diagram
\begin{equation}\label{thecart}
\vcenter{\xymatrix{
T_U\rto^p\dto\ar@{}[dr]|{\Box} & E\dto\\
U\rto^\pi & X}}
\end{equation}
is cartesian,

(v) the composition $T_U\to E\to T_X$ is the derivative of $\pi:U\to X$,

(vi) the  map $E\to \pi\lst T_U$ induced by the cartesian
diagram~(\ref{thecart}), respects brackets (in other  words, $p:T_X\to
E$ is a homomorphism of Lie algebroids covering $\pi:U\to X$).

A local realization is a {\bf realization }(or {\em global
  realization}, if there is fear of confusion) if the structure map
  $\pi:U\to X$ is surjective.
\end{defn}

Note that the existence of global realizations implies surjectivity of
the anchor map.

\begin{rmk}
If we were to ask that $\pi$ be an immersion, instead of a submersion,
we would get essentially the definition of integral submanifold.
\end{rmk}

\begin{defn}\Label{coman}
A {\bf cofoliation }on the manifold $X$ is given by a  Lie algebroid
$a:E\to T_X$ with surjective anchor, which admits a global
realization.  An isomorphic Lie algebroid is considered to define the
same cofoliation.
\end{defn}

\begin{lem}\Label{intlem}
Let $a:E\to T_X$ be a Lie algebroid and let $(V,\rho,r)$ and
$(U,\pi,p)$  be two local realizations. Form the fibred products
$Z=V\times_XU$ and $F=T_V\times_E T_U$. Note that $F\subset T_Z$ is a
subbundle, as $T_Z=T_V\times_{T_X} T_U$.  

Then the  distribution
$F$ on $Z$ is integrable.
\end{lem}
\begin{pf}
Use notation as in the diagram
$$\xymatrix{
Z\dto_q\rto^s\drto^f & U\dto^\pi\\
V\rto^\rho & X}$$
The fact that $F\subset T_Z$ is closed under Lie bracket, can be
checked locally, on a basis of $F$.  Such a local basis may be pulled
back from a local basis for $E$ on $X$, as $F=f\upst E$. Thus, it is
sufficient to prove that $f\lst F\subset f\lst T_Z$ is closed under
the Lie bracket. 

For the following, it is helpful to note that the commutative diagram
of sheaves on $X$
\begin{equation}\Label{help1}
\vcenter{\xymatrix{
f\lst T_Z\rto\dto\ar@{}[dr]|{\Box} & f\lst s\upst T_U\dto\\
f\lst q\upst T_V\rto & f\lst f\upst T_X}}
\end{equation}
is cartesian. Moreover, the diagram
\begin{equation}\Label{help2}
\vcenter{\xymatrix{
f\lst F\rto\dto & f\lst s\upst T_U\dto\\
f\lst q\upst T_V\rto & f\lst f\upst E}}
\end{equation}
is a commutative diagram of isomorphisms.  There is a morphism from
(\ref{help2}) to~(\ref{help1}).

Now let $e_1,e_2$ be two sections of $E$ over $X$. We need to show
that $[f\upst e_1,f\upst e_2]\in f\lst T_Z$ is contained in the subsheaf
$f\lst F$. For this, it suffices to prove that $[f\upst e_1,f\upst
  e_2]$ maps to the same section of $f\lst f\upst E$ under the two
maps $f\lst T_Z\to f\lst s\upst T_U\to f\lst f\upst E$ and $f\lst
T_Z\to f\lst q\upst T_V\to f\lst f\upst E$. In fact, we claim that the
image of $[f\upst e_1,f\upst e_2]$ in $f\lst f\upst E$ is equal to $f\upst
[e_1,e_2]$, under either of these maps. 

Let us prove this for the composition $f\lst T_Z\to f\lst s\upst
T_U\to f\lst f\upst E$.  Note that we have a commutative diagram
$$\xymatrix{
f\lst T_Z\rto & f\lst s\upst T_U\dto & \pi\lst T_U\lto\\
& f\lst f\upst E & E\lto\uto}$$
Now under $E\to\pi\lst T_U$ the section $[e_1,e_2]$ maps to $[\pi\upst
  e_1,\pi\upst e_2]$, because $T_U\to E$ is a morphism of Lie
algebroids. Now we note that $[f\upst e_1,f\upst e_2]$ and $[\pi\upst
  e_1,\pi\upst e_2]$ have same image in $f\lst s\upst T_U$, because
$T_Z\to T_U$ is a morphism of Lie algebroids, and we are done.
\end{pf}

\begin{lem}\Label{banalem}
Let $a:E\to T_X$ be a morphism of vector bundles. 
Let $(U,\pi,p)$ be a triple satisfying items
(i),(ii),\ldots,(v) of Definition~\ref{realize}. Suppose that $\pi$ is
surjective. Form the banal groupoids $U_1=U\times_XU$ and
$E_1=T_U\times_{E}T_U$.  Note that $E_1$ is a connection on $U_1\toto U$. 

The connection $E_1$ is flat if and only if the subsheaf $E\to\pi\lst
T_U$ is closed under the Lie bracket. If this is the case, $E\to T_X$
is a Lie algebroid and $(U,\pi,p)$ a global realization. 
\end{lem}
\begin{pf}
If $E\to\pi\lst T_X$ is closed under the Lie bracket, $E\to T_X$
inherits the structure of a Lie algebroid and $(U,\pi,p)$ is a
realization. From Lemma~\ref{intlem} it follows that $E_1$ is
integrable.  For the converse, note that the diagram 
$$\xymatrix{
f\lst f\upst E & \pi\lst T_U\lto_-{s\upst}\\
\pi\lst T_U\uto^{q\upst} & E\lto\uto}$$
is cartesian, as $\pi$ is surjective. (Here we have borrowed the
notation $s$ and $q$ from the proof of Lemma~\ref{intlem}.) Thus, to
prove that $E$ is closed under the Lie bracket inside $\pi\lst T_U$,
it is enough to prove that for sections $e_1, e_2$ of $E$ the bracket
$[\pi\upst e_1,\pi\upst e_2]$ maps to the same section of $f\lst
f\upst E$ under $s\upst$ and $q\upst$. This is proved by reversing the
argument of Lemma~\ref{intlem}. 
\end{pf}

This Lemma explains what flatness means for connections on banal
groupoids, c.f. Example~\ref{banalex}. 

We also see that a flat connection on a banal groupoid $U_1\toto U_0$
induces a cofoliation on the quotient space $X$. A cofoliation induces
a flat connection on any banal groupoid coming from a local
realization. 

\unsure{
I conjecture that given a pre-realization of a cofoliation, the
difference of the two brackets on $E$ and on $T_{X_0}$ gives rise to a
tensor $B$ inside $\Lambda^2 T_{X_0}\to N=\ker{a}_{X_0}$. The \v Cech
coboundary of $B$ is the curvature of the connection on $X_1$.  One
should probably call $B$ the {\em curving }of the
pre-realization.  A 
realization is a pre-realization whose curving vanishes. Maybe one can
construct the spectral sequence even if the curving is non-zero?
Certainly, it would have to be (covariantly?) closed.  This 
may the be correct integrability condition.

I think it is easy to prove that pre-realizations always exist. Maybe
one can prove that the condition that the curving be closed is
independent of the chosen pre-realization. The closedness of the
curving, might then be used as an axiom to define `integrable
codistribution'. Then one can expect all integrable codistributions
to be cofoliations, in the differentiable and holomorphic context.  In
the algebraic context, there would still be a difference.
}

\subsubsection{Comparing realizations}

Unlike the leaves of a foliation, realizations of
cofoliations are not unique.  Passing to an \'etale cover of $U$, we
get another. We will now develop the substitute for uniqueness.

\begin{defn}\Label{mor-real}
Let $(V,\rho,r)$ and $(U,\pi,p)$  be two local 
realizations for the Lie algebroid $a:E\to T_X$. A {\bf morphism }of
realizations is an \'etale $X$-map $f:V\to U$, such that the induced
diagram
$$\xymatrix@C=1pc{T_V\rrto^{Df}\drto_r && T_U\dlto^p\\
& E &}$$
commutes. 
\end{defn}

The local realizations of a fixed Lie algebroid $E\to T_X$ form thus a
category.  In the differentiable and holomorphic context, this
category is connected:

\begin{lem}\Label{existsthird}
Assume we are in the differentiable or holomorphic category. Let
$a:E\to T_X$ be a Lie algebroid. If $U\to X$ and $V\to X$ are global
realizations of $E$, there exists a third global realization $W\to X$
and morphisms of realizations $W\to U$ and $W\to V$.
\end{lem}
\begin{pf}
Simply form $Z$ as in Lemma~\ref{intlem} and take $W$ to be the
union of enough local leaves of the foliation $F$ to make $W\to X$
surjective. 
\end{pf}

\begin{rmk}
In fact, the category of local realizations behaves a lot like the
category of local \'etale covers of a Deligne-Mumford stack by
manifolds: it admits fibred products and disjoint unions. Any two
objects have a product. (For the product, one has to accept the
disjoint union of the leaves of the foliation $F$ on $Z$ as a
manifold.)
\end{rmk}

\begin{rmk}
One can also take the groupoid with flat connection
$U_1\toto U_0$ induced by any global realization $U_0\to X$ of $E$,
pass to the modified differentiable structure on $U_1$, denoted
$\widetilde{U}_1$, and take the associated `differentiable space'
$\widetilde{X}=[U_0/\widetilde{U}_1]$. If one is willing to accept
$\widetilde{X}$ as some kind of generalized manifold, it will serve
as canonical universal realization.  We could even think of
$\widetilde{X}$ as $X$ with a modified differentiable structure.  This
is in analogy to foliations, where the union of the leaves may be
thought of as a different differentiable structure.  But note that the
dimension of $\tilde{X}$ is equal to $\rank E$, thus {\em larger }than
the dimension of $X$.
\end{rmk}

\subsubsection{Pulling back cofoliations}

Let $a:E\to T_X$ be a cofoliation on the manifold $X$. Let $Y\to X$ be
a smooth map. Define $E_{T_Y}=E\times_{T_X}T_Y$. Note that
$E_{T_Y}$ is  a vector bundle on $Y$, which comes with an epimorphism
of vector bundles $a_Y:E_{T_Y}\to T_Y$. We will define on $E_{T_Y}$ a
Lie algebroid structure making it a cofoliation on $Y$.

For this, choose a global realization $U\to X$ of $E$.
Consider the fibred product $V=U\times_XY$. We have a canonical
morphism $T_V\to E_{T_Y}$, which satisfies Properties
(i),(ii),\ldots,(v) of Definition~\ref{realize}, as can be seen by
contemplating the diagram
$$\xymatrix@=0.5pc{
T_V \ar[rrr] \ar[ddd] \ar[dr] &&& 
E_{T_Y} \ar[dr]\ar[ddd] \ar[rr] && 
T_Y \ar[ddd] \ar [ld] \\ &
V\ar[ddd]\ar[rrr] &&& 
Y \ar[ddd] & \\ &&&&&\\
T_U \ar[dr] \ar[rrr] &&&
E \ar[dr] \ar[rr] &&
T_X\ar[ld] \\ &
U\ar[rrr] &&&
X & }$$
in which the parallelograms with horizontal edges are cartesian. 

Let $E_1$ be the flat connection induced on the banal groupoid
$U_1=U\times_X U$.  Letting $V_1=V\times_Y V$, we have a smooth map
$V_1\to U_1$ and we let $F_1$ be the pullback of the integrable
distribution $E_1$ via $V_1\to U_1$, in other words,
$F_1=E_1\times_{T_{U_1}}T_{V_1}$.  Thus $F_1$ is again an integrable
distribution and, in fact, a flat connection on the banal groupoid
$[V_1\toto V]$. By Lemma~\ref{banalem}, $E_{T_Y}$ is endowed with the
structure of Lie algebroid over $Y$, such that $V\to Y$ is a global
realization. 

\begin{lem}\Label{three-eight}
The Lie algebroid structure on $E_{T_Y}$ is independent of the choice
of the global realization $U\to X$ for $E$.
\end{lem}
\begin{pf}
In the differentiable and holomorphic case, this is easily checked
using Lemma~\ref{existsthird}. In the  algebraic case, one uses the
Lefschetz principle to reduce to the case that the ground field is
$\cc$. Then the equality of two bracket operations can be checked on
the underlying holomorphic manifold. 
\end{pf}

\begin{defn}
We call $a_Y:E_{T_Y}\to T_Y$ with the Lie algebroid structure
constructed above, the
cofoliation on $Y$ obtained by {\bf pull-back }via $Y\to X$. 
\end{defn}

Note that $E_{T_Y}\to E$ is a homomorphism of Lie algebroids covering
$Y\to X$. Of course, pull-back is functorial.

\subsection{Cofoliations on stacks}

We now come to the definition of  cofoliations on stacks. 

Let $X$ be a stack. Let $E$ be a vector bundle over $X$ and $a:E\to
T_X$ an $X$-morphism to the tangent stack. We call $a$ {\em linear},
if for every smooth morphism $U\to X$, where $U$ is a manifold,
forming the fibred product
$$\xymatrix{
E_{T_U}\dto\rto^{a_U}\ar@{}[dr]|{\Box} & T_U\dto\\
E\rto^a\rto & T_X}$$
we obtain a homomorphism of vector bundles $a_U:E_{T_U}\to T_U$. (Note that
$E_{T_U}$ is {\em not }the pullback of $E$ to $U$ via $U\to X$.)

\begin{defn}
A {\bf cofoliation }on the stack $X$ is given by the following data:

(i) a vector bundle $E\to X$,

(ii) a surjective linear map $a:E\to T_X$, (which means that
$a_U:E_{T_U}\to T_U$ is a surjective homomorphism of vector bundles, for
all $U$ as above),

(iii) for every smooth $U\to X$, where $U$ is a manifold, a bracket on
$\Gamma(U,E_{T_U})$. 

\noindent This data is required to satisfy the constraints:

(iv) for every $U$, the bracket makes $E_{T_U}$ a Lie algebroid on $U$,
with anchor map $a_U$,

(v) for any morphism $V\to U$, the induced morphism $E_{T_V}\to
F_{T_U}$ is a homomorphism of Lie algebroids covering $V\to U$,

(vi) the cofoliation is {\em realizable}. By this we mean that there
exists a presentation $\pi:X_0\to X$ of the stack $X$, together
with a morphism of vector bundles $p:T_{X_0}\to E$ covering $\pi$, such
that $p$ identifies $T_{X_0}$ with $\pi\upst E$ and the composition
$a\comp p$ is (isomorphic to!) the canonical morphism $T_{X_0}\to
T_X$. Moreover, we require that  for every smooth $U\to X$, the
fibred product $U_0=X_0\times_X U$ is a realization of the Lie
algebroid $E_{T_U}$. 
\end{defn}

\begin{rmk}
Because cofoliations on manifolds pull back via smooth maps, this
definition of cofoliation on the stack $X$  is equivalent to
Definition~\ref{coman} if $X$ is a manifold. The only reason why we
have to make a new definition for stacks is that we cannot
characterize a cofoliation by the bracket on global sections (or even
\'etale local sections) of $E$.  We need sections over smooth maps
$U\to X$. 
\end{rmk}

\begin{numrmk}\Label{Liebundle}
Note that Condition~(v) implies that the bundle of Lie algebras
$\ker{a_U}$ pulls back to the bundle of Lie algebras $\ker{a_V}$.
Thus we get an induced bundle of Lie algebras $N$ over the stack $X$,
which comes with a homomorphism $N\to E$. We may think of $T_X$ as the
quotient of $E$ by $N$, where $N$ acts on $E$ by addition, through the
homomorphism $N\to E$.  Dually, we have a distinguished triangle
$$\xymatrix@1{L_X\rto &E^\vee \rto & N^\vee\rto & L_X[1]}\,.$$
\end{numrmk}

\begin{rmk}
Let $X$ be a stack and fix a surjective linear map $E\to T_X$, where
$E$ is a vector bundle over $X$. Let $\Xi^p_U$ be the $p$-th exterior
power of the dual vector bundle of $E_{T_U}$.  For fixed $p$, as $U$
varies, the $\Xi^p_U$ form a big sheaf $\Xi^p$ on $X$. Then {\em all }the
compatible brackets on the various $E_{T_U}$ can also be encoded in
{\em one }differential on the big sheaf $\Xi^\ast=\bigoplus_p \Xi^p$.
\end{rmk}

\begin{lem}\Label{third-two}
Definition~\ref{mor-real} and Lemma~\ref{existsthird} carry over to
stacks.
\end{lem}
\begin{pf}
This can be proved by reducing to the manifold case.  The main point
is that one can check the integrability of a foliation smooth
locally. 
\end{pf}

\subsubsection{Relating cofoliations to flat groupoids}

Every flat groupoid defines a cofoliation on the associated stack.

Let $X$ be a stack and fix a surjective linear map $E\to T_X$, where
$E$ is a vector bundle over $X$. We may define a {\bf pre-realization
}of $E\to T_X$ to be a presentation $\pi:X_0\to X$ of the stack $X$
together with $p:T_{X_0}\to E$, such
that $p$ identifies $T_{X_0}$ with $\pi\upst E$ and the composition
$a\comp p$ is the canonical morphism $T_{X_0}\to T_X$.

Any two pre-realizations $X_0\to X$ and $Y_0\to X$ give rise to a
distribution $F\subset T_{Z_0}$, where $Z_0=X_0\times_X Y_0$. In
particular, any pre-realization gives rise to a connection on the
induced groupoid $X_1=X_0\times_X X_0$. 

A pre-realization $X_0\to X$ gives rise to a cofoliation on $X$ if the
induced connection on $X_1$ is integrable. Two such pre-realizations
$X_0\to X$ and $Y_0\to X$ give rise to the same cofoliation on $X$, if
and only if the induced distribution on $Z_0=X_0\times_XY_0$ is
integrable.

In the differentiable or holomorphic context, if the flat groupoids
$X\lcom$ and $Y\lcom$ induce the same cofoliation on $X$, there
exists a third flat groupoid $Z\lcom$ with \'etale horizontal maps
$Z\lcom\to X\lcom$ and $Z\lcom\to Y\lcom$.  The converse is always
true.  

\begin{numex}
Let $\pi:\GG\to X$ be a gerbe over a Deligne-Mumford stack. Then one
may consider cofoliations $a:\pi\upst T_X\to T_{\GG}$ whose anchor
map satisfies $T_\pi\comp a=\id_{T_X}$. (Here $T_{\pi}:T_{\GG}\to T_X$
is the map induced by $\pi$.) These might serve as flat connections on
the `general' gerbe $\GG$. (By `general' we mean not banded by any
group $G$.)  It should be interesting to compare this notion
to the one  defined by Breen and Messing~\cite{BreenMessing}.
\end{numex}

\begin{numrmk}\Label{cryptic}
We did not investigate to what extent realizations always exist. 
It is conceivable, that there exists a `curvature vanishing' condition
that  assures the existence of realizations, at least in the
differentiable and holomorphic categories. 

It is not difficult to prove that pre-realizations always exist. One
might also be able to  construct the spectral
sequence~(\ref{wonder}) using pre-realizations, at least if this
conjectural `curvature  vanishing' condition is satisfied. In that
case, the de Rham differential $d$ would certainly {\em not }square to
zero on the level of co-chains.

This `curvature vanishing' condition might be expressed as the
closedness of a 2-form $B$ as in Remark~\ref{conceivable}.
\end{numrmk}

\subsubsection{Functoriality}

\begin{defn}\Label{maybe-last}
Let $X$ and $Y$ be stacks with cofoliations $E\to T_X$ and $F\to
T_Y$.  A {\bf morphism }of cofoliations is a pair $(f,\phi)$
where $f:X\to Y$ is a morphism of stacks, and $\phi:E\to F$ is a
morphism of vector bundles covering $f$ such that 
$$\xymatrix{
E\rto\dto_{\phi} & T_X\dto^{T_f}\\
F\rto & T_Y}$$
is 2-commutative.  (So $\phi$ involves two implicit 2-morphisms.)
We also ask that there exist realizations $X_0\to X$ of $E$ and
$Y_0\to Y$ of $F$ and a morphism $f_0:X_0\to Y_0$, making 
$$\xymatrix@=0.5pc{
T_{X_0} \ar[rrr] \ar[ddd] \ar[dr] &&& 
E \ar[dr]\ar[ddd] \ar[rr] && 
T_X \ar[ddd] \ar [ld] \\ &
X_0\ar[ddd]\ar[rrr] &&& 
X \ar[ddd] & \\ &&&&&\\
T_{Y_0} \ar[dr] \ar[rrr] &&&
F \ar[dr] \ar[rr] &&
T_Y\ar[ld] \\ &
Y_0\ar[rrr] &&&
Y & }$$
commute, and inducing a horizontal morphism of flat groupoids
$X\lcom\to Y\lcom$.
\end{defn}

\begin{example}
Let $f:X\to Y$ be a representable smooth morphism and $F\to T_Y$ a
cofoliation on $Y$, with realization $Y_0\to Y$. Let $X_0=Y_0\times_Y
X$.  We  get an induced cartesian morphism of groupoids $X\lcom\to
Y\lcom$ and so we can pull back the flat connection on $Y\lcom$ to a
flat connection on $X\lcom$, which defines a cofoliation $E\to T_X$ on
$X$.  By Lemma~\ref{three-eight}, the cofoliation $E$ is
independent of the choice of the presentation $Y_0$ of $F$. 

The  cofoliation $F$ comes with an induced morphism of cofoliations
$E\to F$.  We 
say that $E$ is obtained from $F$ via {\em smooth pullback}. Notation:
$E=f\upst F$. 
\end{example}

\section{The De Rham Complex}

\subsection{The \v Cech complexes of the exterior powers of
  $[\Omega\to\Upsilon]$}

\subsubsection{The trigraded vector space with its first differential $\phi$}

Let $X_1\toto X_0$ be a groupoid with flat connection $E\lcom\to
X\lcom$.  Let $X$ be the associated stack and $N\to E$ the induced
representative of the tangent complex.  Let us denote the dual of
$N\to E$ by 
$$\Omega\stackrel{\phi}{\longrightarrow}\Upsilon\,.$$
It is a homomorphism of vector bundles on $X$.
We also use notation
$\Omega^p=\Lambda^p\Omega$ and $\Upsilon^k=S^k\Upsilon$ (symmetric power). 
We consider the bigraded commutative $\O_{X}$-algebra 
\begin{equation}\Label{bigrad}
\bigoplus_{p\geq
k\geq0}\Omega^{p-k}\otimes\Upsilon^k\,,
\end{equation}
where $\Omega$ is in degree $(1,0)$ and $\Upsilon$ in degree $(1,1)$. We
denote this bigraded algebra by $L=\bigoplus L^{p,k}$, where
$L^{p,k}=\Omega^{p-k}\otimes\Upsilon^k$. 
The homomorphism $\phi$ extends, in a unique way, to a graded
derivation of degree $(0,1)$, which is linear over $\bigoplus_{k\geq0}
\Upsilon^k$.  Explicitly, this graded derivation is given by
$$\phi(\omega_1\wedge\ldots
\wedge\omega_q\otimes\alpha_1\cdot\ldots\cdot\alpha_k) 
= \sum_{i=1}^{q}(-1)^{i+1}\omega_1\wedge\ldots\wedge
\widehat{\omega}_i\wedge 
\ldots\wedge\omega_q\otimes\phi(\omega_i)\cdot
\alpha_1\cdot\ldots\cdot\alpha_k\,,$$
for $\omega_i\in\Omega$ and $\alpha_j\in\Upsilon$. 

Note that $\phi^2=0$, so that $(L,\phi)$ is a differential bigraded
sheaf of $\O_{X}$-algebras.  Via the structure morphism $\pi:X_n\to
X$, we pull back to any $X_n$. 

We now introduce the trigraded $\kk$-vector space
$K=\bigoplus_{p,k,n}K^{p,k,n}$ by
$$K^{p,k,n}=\Gamma(X_n,L^{p,k})=\Gamma(X_n,
\Omega^{p-k}\otimes\Upsilon^k)\,.$$ 
Note that $\phi$ induces a derivation of tridegree $(0,1,0)$ on
$K$ by the formula $\phi=(-1)^n\pi\upst\phi$ and that
$\phi^2=0$. 

\subsubsection{The \v Cech differential}

Since $\hat\pi_q:X_{n+1}\to X_n$ and $\iota_q:X_{n-1}\to X_n$ commute
with the projections to $X$, we  have for each
$q=0,\ldots,n+1$ a homomorphism
\[\hat{\pi}_q\upst:\Gamma({X}_n,\Omega^{p-k}\otimes\Upsilon^k)
\longrightarrow \Gamma({X}_{n+1},\Omega^{p-k}\otimes\Upsilon^k)\]
and each $q=0,\ldots,n-1$ a homomorphism
\[\iota_q\upst:\Gamma({X}_n,\Omega^{p-k}\otimes\Upsilon^k)\longrightarrow
\Gamma({X}_{n-1},\Omega^{p-k}\otimes\Upsilon^k)\,.\]
We can now define the \v Cech differential $\del:K\to K$ of degree
$(0,0,1)$ by
$$\del=\sum_{q=0}^{n+1}(-1)^q\hat{\pi}_q\upst\,.$$
The fact that  $\phi$  is defined over $X$ implies the
following:

\begin{lem}
We have $\del^2=0$, $\phi^2=0$ and $(\del+\phi)^2=0$. \qed
\end{lem}

\subsubsection{The multiplicative structure}

Fix $n,m\geq0$ and consider the two morphisms
\begin{align*}
s:X_{n+m} & \longrightarrow X_n\\
(x_0 \stackrel{\psi_1}{\rightarrow}\ldots
\stackrel{\psi_{n+m}}{\rightarrow} x_{n+m}) & \longmapsto
(x_0 \stackrel{\psi_1}{\rightarrow}\ldots
\stackrel{\psi_{n}}{\rightarrow} x_{n})
\end{align*}
and
\begin{align*}
t:X_{n+m} & \longrightarrow X_m\\
(x_0 \stackrel{\psi_1}{\rightarrow}\ldots
\stackrel{\psi_{n+m}}{\rightarrow} x_{n+m}) & \longmapsto
(x_n \stackrel{\psi_{n+1}}{\rightarrow}\ldots
\stackrel{\psi_{n+m}}{\rightarrow} x_{n+m})\,.
\end{align*}
We define the {\em cup product }on $K$ by
\begin{align*}
\Gamma(X_n,\Omega^{p-k}\otimes\Upsilon^k)\otimes
\Gamma(X_m,\Omega^{p'-k'}\otimes\Upsilon^{k'})& \longrightarrow
\Gamma(X_{n+m},\Omega^{p+p'-k-k'}\otimes\Upsilon^{k+k'})\\
\omega\otimes\omega'&\longmapsto (-1)^{m(p-k)}s\upst(\omega)\wedge
t\upst(\omega')\,. 
\end{align*}
This definition makes $K$ into a trigraded $\kk$-algebra. Both $\phi$ and
$\del$ (and hence also $\del+\phi$) are graded derivations (of total degree
1) with respect to the cup product.

\begin{rmk}
This multiplicative structure is associative. It is not (graded)
commutative. 
\end{rmk}

\subsubsection{The double complex for fixed $p$}

Fixing $p$ and varying $k$ and $n$, we obtain the total complex of a double
complex $(K^{p,\bullet,\bullet},\del,\phi)$.  If $X\lcom$ is
sufficiently nice (all $X_n$ are Stein in the holomorphic context,
affine in the algebraic context), this total complex computes the
hypercohomology of $X\lcom$ with values in $\Lambda^p(\Omega\to\Upsilon)$:
$$h^q(K^{p,\bullet,\bullet},\del+\phi) =
\hh^q\big(X\lcom,\Lambda^p(\Omega\to\Upsilon)\big) =
\hh^q(X,\Lambda^p L_X)\,,$$ 
where $X$ is the stack associated to $X_1\toto X_0$ and $L_X$ its cotangent
complex.

The cup product passes to $\hh(X,\Lambda^pL_X)$, and is commutative on
$\hh(X,\Lambda^pL_X)$.\comment{pf?}

Filtering this double complex by $k$ we obtain
the spectral sequence~(\ref{contemplate}) from the introduction.

\unsure{
\begin{rmk}
The product structure
$$\hh^q\big(X\lcom,\Lambda^p(\Omega\to\Upsilon)\big)\otimes
\hh^{q'}\big(X\lcom,\Lambda^{p'}(\Omega\to\Upsilon)\big)
\longrightarrow
\hh^{q+q'}\big(X\lcom,\Lambda^{p+p'}(\Omega\to\Upsilon)\big)$$
induced by the cup product is equal to the natural product structure
of\comment{reference needed}.
\end{rmk}
}

\subsection{The de Rham differential}

\subsubsection{The alternating de Rham differential $d$}

We  consider an analogue of the bigraded algebra~(\ref{bigrad}), living
on $X_n$. In fact, let us define 
\begin{equation}\Label{second}
L_{X_n}=\bigoplus_{p\geq k\geq 0}L_{X_n}^{p,k}=\bigoplus_{p\geq k\geq
0}\Omega_{X_n}^{p-k}\otimes\Upsilon^k_n\,.
\end{equation}
In keeping with earlier notation, we denote by $L_n$ the pullback of
$L$ to $X_n$. 
We will study how $L_n$ and $L_{X_n}$ relate to each other.  Recall
that, as $n$ varies, all the $L_n$ are pullbacks of each other. In
other words, $L$ is (componentwise) a vector bundle on the stack
$X$. On the other hand, the $L_{X_n}$ do not fit together so nicely.
They only form a big sheaf on $X$.

Recall that $E_n$ is a quotient of $T_{X_n}$ in
$n+1$ different ways, the $q$-th quotient map being denoted
$\eta_q$.  Dually, $\Omega_n$ is a subbundle of
$\Omega_{X_n}$ in $n+1$ different ways.  We denote the embedding dual to
$\eta_q$ again by $\eta_q$. The embedding $E_n\hookrightarrow
T_{X_n}$ corresponds to a quotient map $\delta:\Omega_{X_n}\to
\Omega_n$. Of course, we have that $\delta\comp\eta_q=\id$, for all
$q$. 

Taking alternating powers and tensoring with $\Upsilon_n^k$, we obtain
an embedding $\eta_q:\Omega_n^{p-k}\otimes\Upsilon_n^k\hookrightarrow
\Omega_{X_n}^{p-k}\otimes\Upsilon_n^k$ and a quotient map $\delta$ in the
other direction. We still have $\delta\comp\eta_q=\id$, for all $q$. 
We have thus constructed an algebra morphism $\eta_q:L_n\to L_{X_n}$,
for every $q$, and another algebra morphism $\delta:L_{X_n}\to L_n$,
which is a retraction of every $\eta_q$. 

Now recall that we have the various embeddings $\rho_q:N_n\to
T_{X_n}$.  We denote the dual quotient maps by the same letter:
$\rho_q:\Omega_{X_n}\to\Upsilon_n$. 
The homomorphism $\rho_q$ extends in a unique way to a graded derivation
of degree $(0,1)$ on $L_{X_n}$, which is linear over
$\bigoplus_{k\geq0}\Upsilon_n^k$. We denote this graded derivation
again by $\rho_q$, for all $q=0,\ldots,n$.

\begin{lem}
Letting $\delta_{qr}$ denote the delta of Kronecker, the diagram 
$$\xymatrix{
{\Omega_n^{p-k}\otimes\Upsilon_n^k}\dto_{\eta_r}
\rrto^{\delta_{qr}\phi} && 
{\Omega_n^{p-k-1}\otimes\Upsilon_n^{k+1}}\dto^{\eta_r}\\
{\Omega^{p-k}_{X_n}\otimes\Upsilon_n^k} \rrto^{\rho_{q}} &&
{\Omega^{p-k-1}_{X_n}\otimes\Upsilon_n^{k+1}}}$$
commutes, for all $q,r=0,\ldots,n$. So does the diagram
$$\xymatrix{
{\Omega_n^{p-k}\otimes\Upsilon_n^k} \rrto^{\phi} &&
{\Omega_n^{p-k-1}\otimes\Upsilon_n^{k+1}}\\
{\Omega^{p-k}_{X_n}\otimes\Upsilon_n^k}\uto^{\delta}
\rrto^{\sum_{q=0}^n\rho_{q}} && 
{\phantom{\,.}\Omega^{p-k-1}_{X_n}\otimes\Upsilon_n^{k+1}\,.}
\uto_{\delta} }$$ 
\end{lem}
\begin{pf}
This follows immediately from Diagrams~(\ref{diaqr})
and~(\ref{sumdiagram}).
\end{pf}

Let $\nabla:\Upsilon_0\to\Omega_{X_0}\otimes\Upsilon_0$ be the derived
connection. As remarked, on $\Upsilon_n$, we get $n+1$ different
induced connections 
$\nabla^r:\Upsilon_n\to\Omega_{X_n}\otimes\Upsilon_n$, for
$r=0,\ldots,n$. The $r$-th connection is obtained by thinking of
$\Upsilon_n$ as $\pi_r\upst\Upsilon_0$. To fix notation, Let us focus
on the $0$-th connection $\nabla^0$ and denote it by $\nabla$.  In the
end, our constructions will not depend on this choice. 

We get induced integrable connections on all symmetric powers
$\Upsilon_n^k$ of $\Upsilon_n$:
\[\nabla:\Upsilon_n^k \longrightarrow \Omega_{{X}_n}\otimes 
\Upsilon_n^k\,.\]
Let us denote by 
\[D:\Omega_{X_n}^p\otimes\Upsilon_n^k \longrightarrow
\Omega_{{X}_n}^{p+1}\otimes  \Upsilon_n^k\,\]
the associated covariant derivative, in other words, the differential of
the de~Rham complex of $(\Upsilon_n^k,\nabla)$. 
Note that $D$ is a graded derivation of degree $(1,0)$ on the bigraded
algebra $L_{X_n}$.

\begin{lem}
The derivation $D$ passes to the quotient bigraded algebra $L_n$.
\end{lem}
\begin{pf}
This follows from the fact that $E_n^\perp\subset\Omega_{X_n}$
generates an ideal preserved by the exterior derivative, $E_n$  being
an integrable distribution on $X_n$. 
\end{pf}

We will denote the derivation $D$ induces on $L_n$ by
$$(-1)^nd:\Omega_n^{p-k}\otimes\Upsilon_n^k\longrightarrow
\Omega_n^{p+1-k}\otimes\Upsilon_n^k\,.$$
Taking global sections, we finally arrive at the definition of 
$$(-1)^nd:\Gamma(X_n,\Omega^{p-k}\otimes\Upsilon^k)\longrightarrow
\Gamma(X_n,\Omega^{p+1-k}\otimes\Upsilon^k)\,.$$
Thus, we have defined $d:K^{p,k,n}\to K^{p+1,k,n}$. Note that
$d^2=0$, which follows from $D^2=0$, which holds because $\nabla$ is
an integrable connection.

\begin{lem}\Label{usingother}
Using any of the other connections $\nabla^q$, for $q>0$, gives the
same derivation on the quotient algebra $L_n$. 
\end{lem}
\begin{pf}
This follows immediately from Remark~\ref{afterSR}.
\end{pf}

\begin{cor}
On $\Gamma(X_n,L_n)$, we have $d\del+\del d=0$.
\end{cor}
\begin{pf}
Let us denote the covariant derivative associated with the connection
$\nabla^q$ by $D^q$. Then we have the following relations:
$$\hat\pi_q\upst D^0=\begin{cases} D^1\hat\pi_q\upst&\text{if $q=0$}\,,\\
D^0\hat\pi_q\upst&\text{if $q>0$}\,.\end{cases}$$
The rest is a straightforward calculation, using that $D^0$ and $D^1$
both induce $(-1)^nd$ on the quotient.
\end{pf}

\begin{cor}
The map $d:K\to K$ is a derivation with respect to the cup product of
degree $(1,0,0)$. 
\end{cor}
\begin{pf}
This is a straightforward calculation, but it uses
Lemma~\ref{usingother}.
\end{pf}

\begin{rmk}
Unfortunately, we do not have $d\phi+\phi d=0$. This necessitates the
correction term $\iota$, defined next.
\end{rmk}

\subsubsection{The symmetric de Rham differential, or contraction, $\iota$}

Consider a fixed $X_n$ and define the {\em symmetric partial derivatives}
\[\LL_q:\Omega^p_{{X}_n} \otimes\Upsilon_n^k \longrightarrow
\Omega^{p}_{{X}_n}\otimes\Upsilon_n^{k+1}\,,\]  
for $q=0,\ldots,n$, by 
\[\LL_q=[\rho_q,D]=\rho_q D+D\rho_q\,.\]
Note that the commutator is considered to be a graded commutator, and
both $\rho_q$ and $D$ being odd, we obtain the plus sign in the
formula.  (One may think of $\LL_q$ as the covariant derivative of
$\rho_q$.) 

Define the {\em (total) symmetric derivative }
$$\LL:\Omega^p_{{X}_n} \otimes\Upsilon_n^k \longrightarrow
\Omega^{p}_{{X}_n}\otimes\Upsilon_n^{k+1}\,,\qquad\text{by}\qquad
\LL=\sum_{q=0}^n \LL_q\,.$$
Note that all $\LL_q$, as well as $\LL$, are bigraded derivations of degree
$(0,1)$.

\begin{lem} \Label{sdcs}
The derivations $\LL$ and $\LL_q$, for $q=0,\ldots,n$, pass to the
quotient algebra $\delta:L_{X_n}\to L_n$.
\end{lem}
\begin{pf}
For $\LL$, the claim follows easily from previous results. Therefore,
we only need to prove the claim for $\LL_q$, with $q\geq 1$. 

Let us explain the proof in the differentiable case. 
Recall that we have distinguished the connection
$\nabla=\pi_0\upst\nabla$.  Thus, to do computations, it is convenient
to start from the direct sum decomposition
$$\Omega_{X_n}=\pi_0\upst\Omega_0\oplus
\bigoplus_{j=1}^n\pi_j\upst\Upsilon_0\,.$$ 
This is the dual of the direct sum decomposition~(\ref{dirsu}) for
$q=0$. 
Around a given point of $X_n$ we get an induced decomposition of $X_n$
as a product of $X_0$ with the $\hat\pi_j$-fibres, for
$j=1,\ldots,n$. Since all $\pi_j\upst\Upsilon_0=\Omega_{\hat\pi_j}$ are
endowed with a flat connection, each of these fibres is a flat
manifold.  Choosing arbitrary coordinates $\{y\}$ on $X_0$ and flat
coordinates $\{x_j\}$ on the $\hat\pi_j$-fibre, we get a
coordinate system $\{y,x_1,\ldots,x_n\}$ on $X_n$ with the property
that $dx_j$ is a horizontal frame for $\pi_j\upst\Upsilon_0$, for all
$j=1,\ldots,n$. (Of course, each symbol $x_j$ stands for an
appropriate number of coordinates.)

Now it is easy to describe $\LL_q$, for $q=1,\ldots,n$. We start by
remarking that $\rho_q$ is has been identified with the projection onto
$\pi_q\upst\Upsilon_0$. Since all our coordinate systems $x_j$ are
flat, $\LL_q(dx_j) =0$, whether we think of $dx_j$ as a section of
$\Omega_{X_1}$, or of $\Upsilon_n$. Moreover, $\LL_q(f)=\frac{\partial
  f}{\partial x_q} dx_q$, for a function $f$ on $X_n$. (Now it is also
obvious, why we call the $\LL_q$ partial derivatives.)

To prove the lemma, we need to show that 
$$\LL_q:\Omega_{X_n}\longrightarrow \Omega_{X_n}\otimes \Upsilon_n$$
maps $E_n^\perp\subset \Omega_{X_n}$ into
$E_n^\perp\otimes\Upsilon_n$.  We have identified $E^\perp$ with
$\bigoplus_{j=1}^n\pi_j\upst\Upsilon_0$. 
We have 
$$\LL_q\sum_{j=1}^n f_j\,dx_j=\sum_{j=1}^n dx_j\otimes
\LL_q(f_j)\,,$$
which finishes the proof.  

The proof in the holomorphic case is analogous.  The algebraic case
then follows by appealing to the Lefschetz principle.

For future reference, let us remark, that from our description it
follows that $[\LL_q,\LL_j]=0$, for all $j,q>0$.  
\end{pf}

Because of this lemma it will cause no confusion if we denote the
homomorphisms 
\[\Omega^p_n\otimes\Upsilon_n^k\longrightarrow 
\Omega^p_n\otimes\Upsilon_n^{k+1}\]
induced by $\LL$ and $\LL_q$ by  passing to the quotient, with the
same symbols $\LL$ and $\LL_q$. Moreover, had we used any other
connection $\nabla^r$ to define $\LL_q$, the resulting derivation
would also pass to the quotient. 

\begin{lem}\Label{usingany}
Using any other connection $\nabla^r$, for $r>0$, to define $\LL_q$,
we get the same induced derivation on the quotient $L_n$. 
\end{lem}
\begin{pf}
Recall that $\nabla^r-\nabla^0=\omega_{r0}(\Psi)$. Therefore, we need
to show that $[\omega_{r0}(\Psi),\rho_q]$ kills $\Upsilon_n$ and maps
$\Omega_{X_n}$ into $E_n^\perp\otimes\Upsilon_n$. The latter is clear,
because $\omega_{r0}(\Psi)$ vanishes on $\Omega_{X_n}$ and maps
$\Upsilon_n$ into $E_n^\perp\otimes\Upsilon_n$. 

Now, if $q\not=0,r$, then $\rho_q\comp\omega_{r0}(\Psi)=0$.  If $q=0$
or $q=r$, then $\rho_q\comp\omega_{r0}(\Psi)$ is the symmetrization of
$\Psi$.  But because $\Psi$ is antisymmetric, its symmetrization vanishes.
\end{pf}

\begin{cor} \label{lsps}
For the total symmetric derivative we have 
$$\delta \LL=(\phi d+d\phi)\delta\,.$$
\end{cor}
\begin{pf}
Modulo $E_n^\perp$, we have $\LL=\sum
\LL_q=\sum[\rho_q,D]=[\sum\rho_q,D]=(-1)^n(-1)^n[\phi,d]$. 
\end{pf}

Now consider for each $q=0,\ldots,n-1$ the diagonal
$\iota_q:{X}_{n-1}\rightarrow {X}_n$.  The $p$-the exterior power of the
canonical epimorphism $\iota_q\upst\Omega_{{X}_n}\to\Omega_{{X}_{n-1}}$
gives a map
\[\Gamma({X}_n,\Omega^p_{{X}_n}\otimes\Upsilon_n^k)
\stackrel{{\iota}_q\upst}{\longrightarrow}
\Gamma({X}_{n-1},{\iota}_q\upst\Omega^p_{{X}_{n}}\otimes\Upsilon_{n-1}^k)
\longrightarrow \Gamma({X}_{n-1},\Omega^p_{{X}_{n-1}}\otimes
\Upsilon_{n-1}^k )\,.\]
Note that the diagram
$$\xymatrix{
{\Gamma(X_n,\Omega_{X_n}^{p-k}\otimes\Upsilon^k)}\rrto^{\iota_q\upst}
\dto^{\delta} &&
{\Gamma(X_{n-1},\Omega_{X_{n-1}}^{p-k}\otimes\Upsilon^k)}\dto_{\delta}\\ 
{\Gamma({X}_n,\Omega^{p-k}\otimes\Upsilon_n^k)}
\rrto^{\iota_q\upst} &&
{\Gamma({X}_{n-1},\Omega^{p-k}}\otimes\Upsilon_n^k)}$$
commutes, so there is no ambiguity in the notation $\iota_q\upst$. 

\begin{defn}\Label{contra}
We now define the {\bf contraction}
\[\iota:\Gamma({X}_n,\Omega^{p-k}\otimes\Upsilon^k) \longrightarrow
\Gamma({X}_{n-1},\Omega^{p-k}\otimes\Upsilon^{k+1})\]
by 
\[-\iota=\sum_{0\leq i< j\leq n}(-1)^i\iota_i\upst \LL_j\,.\]
\end{defn}

Note that $\iota:K^{p,k,n}\to K^{p+1,k+1,n-1}$ is of degree $(1,1,-1)$. 

\begin{prop}
$[d,\iota]=0$.
\end{prop}
\begin{pf}
A simple calculation shows that $[D,\LL_j]=0$, for all $j$. Since the
derivative $D$ commutes with restriction via $\iota_j\upst$, it
follows that $D\iota=\iota D$.  The formula of the proposition follows
by passing to the quotient.
\end{pf}

\begin{prop}
$[\phi,\iota]=0$.
\end{prop}
\begin{pf}
It suffices to  prove that $[\phi,\LL_j]=0$, for all $j>0$. Clearly,
$[\phi,\LL_j]$ vanishes 
on $\bigoplus\Upsilon_n^k$. So it is enough to check vanishing of
$[\phi,\LL_j]$ on a
set of generators of $\Omega_{X_n}$ over $\O$. Since only $\LL_j$ for
$j>0$ are involved, we may use coordinates as in the proof of
Lemma~\ref{sdcs}.  Because we are proving an identity on the quotient, it
suffices to check the identity on the generators $dy$. Because all $j$
involved are greater than zero, we know that $\LL_j(dy)=0$ and hence
$\phi\LL_j(dy)=0$. Also note that $\phi(dy)=\pi_0\upst\phi(dy)$, and
so $\phi(dy)$ is constant in the $x_j$-directions.  Hence
$\LL_j\phi(dy)=0$. Adding up, we get the required formula
$[\phi,\LL_j](dy)=0$. 
\end{pf}

\begin{lem}\Label{rhoi}
For all $q=0,\ldots,n-1$ and all $j=0,\ldots,n-1$, we have
$$\rho_q\iota_j\upst=\begin{cases}
\iota_j\upst\rho_q &\text{if $q<j$}\,,\\
\iota_j\upst(\rho_{q}+\rho_{q+1})&\text{if $q=j$}\,,\\
\iota_j\upst\rho_{q+1} & \text{if $q>j$}\,.\end{cases}$$
\end{lem}
\begin{pf}
Let us first consider the cases $q\not=j$. Then we have the cartesian
diagram
$$\xymatrix{
X_n\rto^{\hat\pi_{q'}} & X_{n-1}\\
X_{n-1}\uto^{\iota_j}\rto^{\hat\pi_{q}} & X_{n-2}\uto_{\iota_{j'}}}$$
If $q<j$, then $q'=q$ and $j'=j-1$. If $q>j$, then $q'=q+1$ and
$j'=j$. It shows that we have a commutative diagram
$$\xymatrix{
\iota_j\upst\Omega_{X_n}\dto\rto^{\iota_j\upst\rho_{q'}}&
\iota_j\upst\Omega_{\hat\pi_{q'}}\dto^=\\
\Omega_{X_{n-1}}\rto^{\rho_q}& \Omega_{\hat\pi_{q}}}$$
which proves the formula for the case $q\not=j$.  For the remaining
case, notice that $\sum\rho_q=\phi\delta$, so that $\sum\rho_q$ commutes
with $\iota_j$ and we have
$(\sum_{q=0}^{n-1}\rho_q)\iota_j\upst=\iota_j\upst\sum_{q=0}^n\rho_q$.
Combine this with the case $q\not=j$ to finish.
\end{pf}

\begin{cor}\Label{rhol}
For all $q=0,\ldots,n-1$ and all $j=0,\ldots,n-1$, we have
$$\LL_q\iota_j\upst=\begin{cases}
\iota_j\upst\LL_q &\text{if $q<j$}\,,\\
\iota_j\upst(\LL_{q}+\LL_{q+1})&\text{if $q=j$}\,,\\
\iota_j\upst\LL_{q+1} & \text{if $q>j$}\,.\end{cases}$$
\end{cor}

\begin{prop} \Label{bipd}
$\iota^2=0$.
\end{prop}
\begin{pf}
This is a straightforward calculation using Corollary~\ref{rhol} and
the fact that $[\LL_q,\LL_j]=0$ on $X_n$, if $n\geq2$, which we
remarked in the proof of Lemma~\ref{sdcs}.
\end{pf}

\begin{lem}
Consider $\hat\pi_q:X_{n+1}\to X_n$. Then for every $j=0,\ldots,n+1$,
we have
$$\rho_j\hat\pi_q\upst=\begin{cases}
\hat\pi_q\upst\rho_j & \text{if $j<q$}\,,\\
0 & \text{if $j=q$}\,,\\
\hat\pi_q\upst\rho_{j-1} & \text{if $j>q$}\,.\end{cases}$$
\end{lem}
\begin{pf}
The proof is analogous to the proof of Lemma~\ref{rhoi}, including the
trick to reduce to the case $j\not=q$.
\end{pf}

\begin{cor}\Label{lastco}
We have
$$\LL_j\hat\pi_q\upst=\begin{cases}
\hat\pi_q\upst\LL_j & \text{if $j<q$}\,,\\
0 & \text{if $j=q$}\,,\\
\hat\pi_q\upst\LL_{j-1} & \text{if $j>q$}\,.\end{cases}$$
\end{cor}
\begin{pf}
Comparing with the proof of Corollary~\ref{rhol}, there is the added
subtlety that $\hat\pi_0\upst$ does not commute with $D$. So, at least
on $L_{X_n}$, we do not have
$\LL_1\hat\pi_0\upst=\hat\pi_0\upst\LL_0$. But because of
Lemma~\ref{usingany}, if we pass to the quotient $L_n$, we do have
this required equality.
\end{pf}

\begin{prop} \label{cdiy}
$[\del,\iota]=-\LL$.
\end{prop}
\begin{pf}
This is a straightforward calculation using Corollary~\ref{lastco}.
\end{pf}

\begin{rmk}
Unfortunately, $\iota:K\to K$ is not a
derivation with respect to the cup product.  Rather, we have the
formula
$$\iota(\alpha\cup\beta)=\iota\alpha\cup\beta+(-1)^{\deg
  \alpha}\alpha\cup \iota\beta - I(\alpha)\cup \LL\beta\,.$$
Here $I=\sum_i(-1)^i\iota_i\upst$. This follows easily from repeated
  applications of Corollary~\ref{lastco}.
The error term can be expressed as follows:
\begin{multline}\Label{error}
-I\alpha\cup\LL\beta\\
=I\alpha\cup\iota(\phi+\del)\beta
-(-1)^{\deg\alpha} I(\phi+\del)\alpha\cup\iota\beta
-(-1)^{\deg\alpha} (\phi+\del)(I\alpha\cup\iota\beta)\,.
\end{multline}
Thus, after passing to cohomology with respect to $(\phi+\del)$, this
error term vanishes.
\end{rmk}

\begin{cor}
Let us summarize.  We have four maps on the complex $K$ for which
the following formulae hold: 
\begin{enumerate}
\item $\phi^2=\del^2=d^2=\iota^2=0$.
\item $[\phi,\del]=[\phi,\iota]=[\del,d]=[d,\iota]=0$.
\item $[\phi,d]+[\del,\iota]=0$.
\item $(\phi+\del)^2=(d+\iota)^2=[\phi+\del,d+\iota]=0$.
\item $(\phi+\del+d+\iota)^2=0$.
\end{enumerate}
The degrees are as follows:
$$\deg\phi=(0,1,0)\quad\deg\del=(0,0,1)\quad\deg d=(1,0,0)\quad
\deg\iota=(1,1,-1)\,.$$
Finally, $\phi$, $\del$ and $d$ are derivations, whereas $\iota$ is a
derivation up to the error term~(\ref{error}). 

We will use $K$ to denote this trigraded $\kk$-vector space with the
differential $(\phi+\del+d+\iota)$. If we need to include
the flat groupoid  $(X\lcom,E\lcom)$ in the notation,
we denote this complex by $K(X\lcom,E\lcom)$, or simply $K(X\lcom)$.
\end{cor}

\unsure{

\begin{lem} \label{iaph}
Let $0\leq q\leq p+1$ and $0\leq r\leq p$. Then for endomorphisms of
$\Gamma({X}_p,\Omega_p^n\otimes \Upsilon^k_p)$ the following
relations hold. 

(i) If $q<r$, then
$\iota_r\upst\hat{\pi}_q\upst=\hat{\pi}_q\upst\iota_{r-1}\upst$. 

(ii) If $q>r+1$, then
$\iota_r\upst\hat{\pi}_q\upst=\hat{\pi}_{q-1}\upst\iota_{r}\upst$. 

(iii) If $q=r$ or $q=r+1$, then $\iota_r\upst\hat{\pi}_q\upst=\id$.
\end{lem}
\begin{pf}
This follows from corresponding relations for the morphisms $\iota_r$
and $\hat{\pi}_q$ themselves. 
\end{pf}

}

\begin{prop}[naturality]\Label{pre-natural}
Given a horizontal morphism of flat groupoids 
$f:X\lcom\to Y\lcom$, we get an induced morphism $f\upst:K(Y\lcom)\to
K(X\lcom)$. This morphism respects the triple grading, the cup product
and all four differentials $\phi$, $\del$, $d$ and $\iota$. \qed
\end{prop}

\subsection{Conclusions}

\begin{them}
Let $X$ be a stack.  If we are in the holomorphic category, assume
that the diagonal of $X$ is relatively Stein, if we are in the
algebraic category, assume that $X$ has affine diagonal. 

Any cofoliation on $X$ gives rise to
an $E_1$-spectral sequence 
\begin{equation}\label{not-again}
E_1^{m,n}=\bigoplus_{p+k=m}
H^n(X,\Omega^{p-k}\otimes\Upsilon^k)\Longrightarrow
H^{m+n}_{DR}(X)\,.
\end{equation}
For every $n$, the term $E^{\bullet,n}_1$ is a double complex, whose
differentials are  induced by 
$\phi$ and $d$, respectively.
\end{them}
\begin{caveat}[algebraic case]
The techniques developed in this paper do not suffice to prove in the
algebraic case that the spectral sequence~(\ref{not-again}) is
an invariant of the cofoliation.  Thus, in the algebraic case, it
remains a conjecture that 
this spectral sequence does not
depend on a chosen flat groupoid realization.
\end{caveat}
\begin{pf}
By our
assumptions on the stack $X$, we can choose a flat groupoid $X\lcom$
realizing our cofoliation on $X$, 
where every $X_n$ is a differentiable manifold, a Stein holomorphic
manifold or a smooth affine variety, depending on the context.
The purpose of this is so that
we can use \v Cech cohomology of the simplicial manifold $X\lcom$ to
compute cohomology of $X$ with values in coherent sheaves over $X$.

Now we filter the complex
$K(X\lcom)$ by the degree $m=p+k$. This gives rise to a spectral
sequence abutting to the (total) cohomology of $K(X\lcom)$, with the
given $E_1$-term. 
We need to prove two things: 

(i) the cohomology of $K(X\lcom)$ is equal to the de Rham cohomology
of $X$, 

(ii) the spectral sequence is independent of the flat groupoid
$X\lcom$ chosen to realize the given cofoliation on $X$.

Let us prove (i).
Note that the flat connection on $X_1\toto X_0$ induces flat
connections on every `shifted' groupoid $X_{n+1}\toto X_n$, by {\em
  cartesian pullback}.  If we 
choose $X_{n+1}\toto X_n$ to have source and target maps $\hat\pi_q$
and $\hat\pi_r$, then this flat connection is given by
$E_{n+1}\oplus\bigoplus_{j\not=q,r}\pi_j\upst 
N\subset T_{X_{n+1}}$.

\newcommand{\rtto}{\ar@<-.5ex>[r]\ar@<.5ex>[r]}
\newcommand{\rttto}{\ar[r]\ar@<1ex>[r]\ar@<-1ex>[r]}
\newcommand{\rtttto}{\ar@<-1.5ex>[r]\ar@<-.5ex>[r]\ar@<1.5ex>[r]\ar@<.5ex>[r]}
\newcommand{\dtto}{\ar@<-.5ex>[d]\ar@<.5ex>[d]}
\newcommand{\dttto}{\ar[d]\ar@<1ex>[d]\ar@<-1ex>[d]}
\newcommand{\dtttto}{\ar@<-1.5ex>[d]\ar@<-.5ex>[d]\ar@<1.5ex>[d]\ar@<.5ex>[d]} 

Let us fix notation as follows: $X\lcom[n-1]=[X_{n+1}\toto X_n]$, with
source and target maps $\hat\pi_0$ and $\hat\pi_1$.  Then we have a
diagram
$$\xymatrix@1{
\ldots\rtttto & X\lcom[2]\rttto & X\lcom[1] \rtto & X\lcom[0]\rto &
X\lcom}$$
Of horizontal morphisms of flat groupoids. We get
an augmented cosimplicial diagram
$$\xymatrix@1{
K(X\lcom)\rto & K(X\lcom[0])\rtto & K(X\lcom[1])\rttto &
K(X\lcom[2])\rtttto & \ldots}$$
of morphisms of complexes. It is not difficult to check that if we
pass to the alternating sums of the horizontal maps, we get a
resolution of $K(X\lcom)$:
$$\xymatrix@1{
K(X\lcom)\rto & K(X\lcom[0])\rto^{\del} & K(X\lcom[1])\rto^\del &
K(X\lcom[2])\rto^\del & \ldots}$$
Thus, if we assemble all $K(X\lcom[\ast])$ into a
double complex with the \v Cech differential $\del$ we end up with a
quasi-isomorphism of complexes
\begin{equation}\Label{firstclaim}
K(X\lcom)\longrightarrow K(X\lcom[\ast])\,.
\end{equation}
Now each $X\lcom[n]$ is the banal groupoid associated to a surjective
submersion $X_{n+1}\to X_n$.  Thus we have a horizontal morphism
$X\lcom[n]\to X_n$, where $X_n$ stands for the groupoid $X_n\toto
X_n$ with its canonical flat connection. We get an induced
quasi-isomorphism $K(X_n)\to K(X\lcom[n])$, as 
both groupoids present the manifold $X_n$.  But $K(X_n)$ is the usual
de Rham complex of $X_n$, so that we have, in fact, a
quasi-isomorphism $\Omega\com(X_n)\to K(X\lcom[n])$.  Assembling all
these quasi-isomorphism together, we obtain another quasi-isomorphism
$$\Omega\com(X\lcom)\longrightarrow K(X\lcom[\ast])\,.$$
Together with (\ref{firstclaim}), we see that the cohomology of
$\Omega\com(X\lcom)$ and $K(X\lcom)$ are canonically isomorphic.
Since $\Omega\com(X\lcom)$ computes the de Rham cohomology of $X$,
this completes the proof of~(i). 

To prove (ii), we start by noticing that an \'etale morphism of
groupoids realizing the given cofoliation induces the identity on the
$E_1$-term of our spectral sequence.  Hence any two morphisms between
the same two realizations induce the same canonical isomorphism of our
spectral sequence. This proves~(ii) in the differentiable and
holomorphic case, by appealing to Lemma~\ref{third-two}.  In the
algebraic case, one might prove independence of the realization by
following the program outlined in Remark~\ref{cryptic}.
\end{pf}

\begin{prop}[naturality]
The spectral sequence (\ref{not-again}) is natural for morphisms of
cofoliations (see Definition~\ref{maybe-last}). 
In particular, it commutes with smooth pullbacks, hence with \'etale
localization. 
\end{prop}
\begin{pf}
Since every morphism of cofoliations is induced by a flat morphism of
groupoids, we can use Proposition~\ref{pre-natural}, to obtain an
induced morphism of spectral sequences. We need to prove that this
this morphism of spectral sequences is independent of the realizations
chosen. This is enough to check on the $E_1$-level.  But on the
$E_1$-level, any flat morphism of realizations induces the
homomorphism $H^n(Y,\Lambda^{p-k}F^\vee\otimes S^k M^\vee)\to
H^n(X,\Lambda^{p-k}E^\vee\otimes S^k N^\vee)$ coming from
$(f,\phi):E\to F$. (Here $M=\ker(F\to T_Y)$ and 
$N=\ker(E\to T_X)$. The  other notation is taken from
Definition~\ref{maybe-last}.)
\end{pf}

\begin{prop}[multiplicativity]
The spectral sequence~(\ref{not-again}) is multiplicative, i.e.,
consists of 
differential graded $\kk$-algebras and derivations.
\end{prop}
\begin{pf}
Let $\widetilde{K}\subset K$ be the normalization of $K$. This is the
subcomplex of elements vanishing under all pullback maps
$\iota_j$. Nothing changes, except now $\iota$ and hence
$(\phi+\del+d+\iota)$ is a derivation.  The proposition follows.
\end{pf}

\begin{rmk}
The spectral sequence~(\ref{not-again}) has more structure.  For
example, we can filter the  whole spectral sequence using the degree
$p$. In other words, we construct a spectral sequence of filtered
algebras. Since the $E_\nu$-level is now a filtered complex, it has an
associated spectral sequence.  It's abutment is the $E_\nu$-term of
the Hodge to De Rham spectral sequence~(\ref{it}). 
\end{rmk}


\end{document}